\documentclass[letterpaper, 11pt]{article}
\usepackage[utf8]{inputenc}

\usepackage{geometry}
\usepackage{amssymb}
\usepackage{amsmath}
\usepackage{breqn}
\usepackage{enumerate}
\usepackage{algorithm}
\usepackage{algpseudocode}
\usepackage{amsthm}
\usepackage{hyperref}
\usepackage{tikz}
\usetikzlibrary{shapes.geometric, arrows, backgrounds, scopes, positioning}
\usepackage{graphicx}
\usepackage{subcaption}
\usepackage{colortbl}
\usepackage{cite}
\usepackage{authblk}

\def\R{\mathbb{R}}
\def\A{A^TA}

\def\Ai{\left(a^{i}\right)^T a^{i}}
\def\Azi{\left(a^{\zeta_{t_t}}\right)^T a^{\zeta_{t_t}}}
\def\E{\mathbb{E}}
\def\Et{\mathbb{E}_{t}}
\def\Eit{\mathbb{E}_{i_t}}
\def\Ezt{\mathbb{E}_{\zeta_t}}

\def\V{\mathbb{V}_{i_t}}
\def\Ax{\mathcal{A}}

\providecommand{\norm}[1]{\ensuremath{\left\lVert#1\right\rVert }}
\providecommand{\mnorm}[1]{\ensuremath{\left\lvert#1\right\rvert}}
\providecommand{\condexp}[1]{\Eit\ensuremath{\left[#1\right]}}
\providecommand{\serverexp}[1]{\Ezt\ensuremath{\left[#1\right]}}
\providecommand{\agentexp}[1]{\mathbb{E}_{I_t}\ensuremath{\left[#1\right]}}
\providecommand{\totexp}[1]{\Et\ensuremath{\left[#1\right]}}

\providecommand{\var}[1]{\V\ensuremath{\left[#1\right]}}
\DeclareRobustCommand{\bigO}{%
  \text{\usefont{OMS}{cmsy}{m}{n}O}%
}

\newtheorem{theorem}{\bfseries Theorem}
\newtheorem*{theorem*}{\bfseries Theorem}
\newtheorem{lemma}{\bfseries Lemma}
\newtheorem{assumption}{\bfseries Assumption}

\title{
Accelerating Distributed SGD for Linear Regression \\using Iterative Pre-Conditioning}

\author{Kushal Chakrabarti$^\star$, Nirupam Gupta$^\dagger$, and Nikhil Chopra$^\star$
\thanks{$^\star$ University of Maryland, College Park, Maryland 20742, U.S.A. \\
$^\dagger$ Georgetown University, Washington, DC 20057, U.S.A. \\
Emails: {\em kchak@terpmail.umd.edu}, {\em nirupam115@gmail.com}, and {\em nchopra@umd.edu}}%
}

\date{}

\begin{document}

\maketitle

\begin{abstract}

This paper considers the multi-agent distributed linear least-squares problem. The system comprises multiple agents, each agent with a locally observed set of data points, and a common server with whom the agents can interact. The agents' goal is to compute a linear model that best fits the collective data points observed by all the agents. In the server-based distributed settings, the server cannot access the data points held by the agents. The recently proposed Iteratively Pre-conditioned Gradient-descent (IPG) method~\cite{chak2020ipc} has been shown to converge faster than other existing distributed algorithms that solve this problem. In the IPG algorithm, the server and the agents perform numerous iterative computations. Each of these iterations relies on the entire batch of data points observed by the agents for updating the current estimate of the solution. Here, we extend the idea of iterative pre-conditioning to the {\em stochastic} settings, where the server updates the estimate and the {\em iterative pre-conditioning matrix} based on a single randomly selected data point at every iteration. We show that our proposed Iteratively Pre-conditioned Stochastic Gradient-descent (IPSG) method converges linearly in expectation to a proximity of the solution. Importantly, we empirically show that the proposed IPSG method's convergence rate compares favorably to prominent stochastic algorithms for solving the linear least-squares problem in server-based networks.

\end{abstract}

\section{INTRODUCTION}
\label{sec:intro}

\tikzstyle{master} = [rectangle, rounded corners, minimum width=1.7cm, minimum height=1cm,text centered, text width=1cm, draw=black, fill=blue!30]
\tikzstyle{machine} = [rectangle, minimum width=1cm, minimum height=1cm,text centered, text width=2cm, draw=black, fill=blue!10]
\tikzstyle{dots} = [circle, inner sep=0pt,minimum size=2pt, draw=black, fill=blue!50!cyan]
\tikzstyle{arrow} = [thick,<->,>=stealth]

\begin{figure}[thpb]  
\centering
\begin{tikzpicture}[node distance = 1.5cm, auto]
    \node (master) [master] {Server};
    \node (m/c2) [machine, below  = of master] {Agent 2\\ $(A^2,B^2)$};
    \node (m/c1) [machine, left = of m/c2, xshift=1cm] {Agent 1\\ $(A^1,B^1)$};
    \node (d1) [dots, right = of m/c2, xshift=-0.8cm] {};
    \node (d2) [dots, right = of d1, xshift=-1.4cm] {};
    \node (d3) [dots, right = of d2, xshift=-1.4cm] {};
    \node (m/c3) [machine, right = of d3, xshift=-0.8cm] {Agent m\\ $(A^m,B^m)$};
    
    \draw[arrow] (master) -- (m/c1);
    \draw[arrow] (master) -- (m/c2);
    \draw[arrow] (master) -- (m/c3);
\end{tikzpicture}
\caption{System architecture.}
\label{fig:sys}
\end{figure}
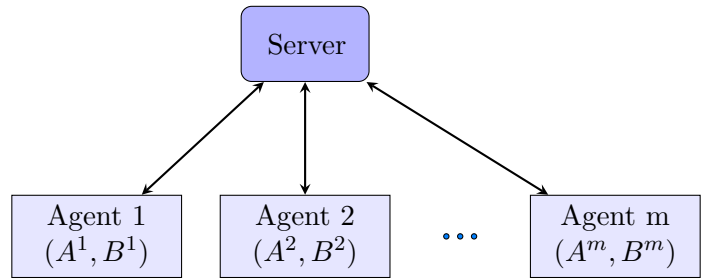

This paper considers solving the distributed linear least-squares problem using stochastic algorithms. In particular, as shown in in Fig.~\ref{fig:sys}, we consider a server-based distributed system comprising $m$ agents and a central server. The agents can only interact with the server, and the overall system is assumed synchronous. Each agent $i$ has $n$ {\em local} data points, represented by an {\em input matrix} $A^i$ and an {\em output vector} $B^i$ of dimensions $n \times d$ and $n \times 1$, respectively. Thus, for all $i \in \{1, \ldots, \, m\}$, $A^i \in \R^{n \times d}$ and $B^{i} \in \R^{n}$. 
For each agent $i$, we define a {\em local cost function} $F_i: \R^d \to \R$ such that for a given {\em regression parameter} vector $x \in \R^d$,
\begin{align}
     F^i(x) = \frac{1}{2n} \norm{A^i x - B^i}^2, \label{eqn:dist}
\end{align}
where $\norm{\cdot}$ denotes the Euclidean norm. The agents' objective is to compute an optimal parameter vector $x^* \in \R^d$ such that
\begin{align}
     x^* \in  \arg \min_{x \in \R^d} \, \frac{1}{m}\sum_{i=1}^m F^i(x). \label{eqn:opt_1}
\end{align}
Since each agent knows only a segment of the {\em collective data points}, they collaborate with the server for solving the distributed problem~\eqref{eqn:opt_1}. However, the agents \underline{do not} share their local data points with the server. An algorithm that enables the agents to jointly solve the above problem in the architecture of Fig.~\ref{fig:sys} without sharing their data points is defined as a {\em distributed algorithm}. \\

%The above distributed linear least-squares problem finds its application in a wide range of real-world problems, including state estimation in sensor networks~\cite{}, training supervised machine learning models~\cite{}, virtual service generation in fault tolerant internet of things (IoT)~\cite{}, quality of experience evaluation in IoT~\cite{}, computed tomography reconstruction~\cite{}. 
%The data points for many of these applications exist as dispersed over a number of isolated sources (agents). This intrinsic distributed nature of the data points by default requires such computational task to be solved in distributed manner. Even when the entire dataset is available to a single agent, large number of data points makes it challenging to solve the problem on a single machine because of its limited memory and computational capability. Such large-scale tasks usually rely on distributed or parallel computations over multiple processing machines (agents), such as~\eqref{eqn:opt_1}. 

There are several theoretical and practical reasons for solving the distributed problem~\eqref{eqn:opt_1} using stochastic methods rather than batched optimization methods, particularly when the number of data-points is abundant~\cite{bottou2018optimization}. The basic prototype of the stochastic optimization methods that solve~\eqref{eqn:opt_1} is the traditional stochastic gradient (SGD)~\cite{bottou2018optimization}. Several accelerated variants of the stochastic gradient descent algorithm have been proposed in the past decade~\cite{duchi2011adaptive, kingma2014adam, zeiler2012adadelta,tieleman2012lecture, reddi2019convergence, dozat2016incorporating}. A few of such well-known methods are the  adaptive gradient descent (AdaGrad)~\cite{duchi2011adaptive}, adaptive momentum estimation (Adam)~\cite{kingma2014adam}, AMSGrad~\cite{reddi2019convergence}. These algorithms are stochastic, wherein the server maintains an estimate of a solution defined by~\eqref{eqn:opt_1}, which is refined iteratively by the server using the {\em stochastic gradients} computed by a randomly chosen agent. \\

In particular, Adam has been demonstrated to compare favorably with other stochastic optimization algorithms for a wide range of optimization problems~\cite{radford2015unsupervised, peters2018deep, wu2016google}. However, Adam updates the current estimate effectively based on only a window of the past gradients due to the exponentially decaying term present in its estimate updating equation, which leads to poor convergence in many problems~\cite{reddi2019convergence}. 
A recently proposed variant of Adam is the AMSGrad algorithm, which proposes to fix Adam's convergence issue by incorporating ``long-term memory'' of the past gradients. \\
% While AMSGrad has been shown to perform better than Adam on CIFAR-10 dataset~\cite{reddi2019convergence}, other experiments suggests AMSGrad to be similar or worse than Adam~\cite{}. 

In this paper, we propose a {\em stochastic iterative pre-conditioning} technique for improving the rate of convergence of the distributed stochastic gradient descent method when solving the linear least-squares problem~\eqref{eqn:opt_1} in distributed networks. The idea of {\em iterative pre-conditioning} in the deterministic (batched data) case has been first proposed in~\cite{chak2020ipc}, wherein the server updates the estimate using the sum of the agents' gradient multiplied with a suitable iterative pre-conditioning matrix. Updating the pre-conditioning matrix depends on the entire dataset at each iteration.
The proposed algorithm extends that idea to the {\em stochastic} settings, where the server updates both the estimate and the iterative pre-conditioning matrix based on a randomly chosen agents' {\em stochastic gradient} at every iteration. Each agent computes its stochastic gradient based on a single randomly chosen data point from its local set of data points. Using real-world datasets, we empirically show that the proposed algorithm \underline{converges in fewer iterations} compared to the aforementioned state-of-the-art distributed methods. \\

We note that the prior work on the formal convergence of the iteratively pre-conditioned gradient-descent method only considers the batched data at every iteration~\cite{chak2020ipc}. Therefore, besides empirical results, we also present a formal analysis of the proposed algorithm's convergence in stochastic settings. \\

Details of our proposed algorithm are presented in Section~\ref{sec:algo}, and the formal results on its convergence are in Section~\ref{sec:conv}. Our empirical results, including comparisons of our proposed algorithm with the related state-of-the-art distributed stochastic methods, are presented in Section~\ref{sec:exp}. Below, we present a summary of our key contributions.

% \subsection{Extension to Distributed Logistic Regression}
% \label{sub:logistic_def}

% We have considered the multi-agent distributed linear least-squares problem~\eqref{eqn:opt_1} where each agent's local cost function $F^i(x)$ is in squared form and they collaborate with a server to minimize the aggregate cost function $\sum_{i=1}^m F^i(x)$. Below we consider the case where the local cost function $F^i$ of each agent $i \in \{1,\ldots,m\}$ is defined as the log-likelihood of the {\em multinomial logit function}~\cite{}. Specifically, for a given parameter matrix $X = \begin{bmatrix} x_1, \ldots, \, x_m \end{bmatrix} \in \R^{d \times K}$ with each column $j\in \{1,\ldots,m\}$ of $X$ denoted as $x_j \in \R^d$,
% \begin{align}
%     F^i(X) & = - \sum_{j=1}^n \sum_{k=1}^K \textbf{1}\{B^i_j=k\} \log \frac{\exp\left((A^i_j)^T x_k \right)}{\sum_{l=1}^K \exp\left((A^i_j)^T x_j \right)}. \label{eqn:dist_logit}
% \end{align}
% Here, $\textbf{1}\{\cdot\}$ denotes the indicator function and $K$ is the number of classes. The agents' objective is to compute an optimal parameter matrix $X^* \in \R^{d \times K}$ such that
% \begin{align}
%     X^* \in  \arg \min_{X \in \R^{d \times K}} \, \frac{1}{N}\sum_{i=1}^N F^i(X). \label{eqn:opt_logit}
% \end{align}

% The proposed algorithm can be extended to solve the above optimization problem. For experimental evidence on real-world datasets, please refer to Section~\ref{}.

\subsection{Summary of Our Contributions}
\label{sub:contri}

\begin{enumerate}[(i)]
    \item We present a formal convergence analysis of our proposed algorithm. Our convergence result can be informally summarized as follows. Formal details are presented in Theorem~\ref{thm:z_conv} in Section~\ref{sub:theorem}. 
    % in our key result, Theorem~\ref{thm:z_conv}, which is informally stated below. 
    \begin{theorem*}[{\it informal}]
    % Consider Algorithm~\ref{algo_1} with positive parameters $\alpha$, $\delta$, and $\beta$. 
    Suppose the solution of problem~\eqref{eqn:opt_1} is unique, and the variances of the stochastic gradients computed by the agents are bounded. In that case, our proposed algorithm, i.e., Algorithm~\ref{algo_1}, converges {\em linearly} in expectation to a proximity of the solution of the problem~\eqref{eqn:opt_1}.  
    % If the parameters $\alpha$ and the stepsize $\delta$ are sufficiently small, then for arbitrary initialization, Algorithm~\ref{algo_1} . 
    % bound of the gradients' variance, as well as the parameter  and .
    \end{theorem*}
    
    The approximation error is proportional to the parameter values $\alpha$, $\delta$, and the variances of the stochastic gradients. It should be noted that, as shown in prior work~\cite{chakrabarti2020iterative}, our algorithm converges {\em superlinearly} to the exact solution when the gradient noise is zero. 
    % at the minimum point, our algorithm converges linearly to the exact minimum point, as expected for a stochastic algorithm.  
     
    \item Using real-world datasets, we empirically show that our proposed algorithm's convergence rate is superior to that of the state-of-the-art stochastic methods when distributively solving linear least-squares problems. These datasets comprise
    \begin{itemize}
        \item four benchmark datasets from the SuiteSparse Matrix Collection;
        \item a subset of the \textit{``cleveland''} dataset from the UCI Machine Learning Repository, which contains binary classification data of whether the patient has heart failure or not based on 13 features;
        \item a subset of the \textit{``MNIST''} dataset for classification of handwritten digit one and digit five.
    \end{itemize}
    Please refer to Section~\ref{sec:exp} for further details.
\end{enumerate}

\section{SGD WITH ITERATIVE PRE-CONDITIONING}
\label{sec:algo}

In this section, we present our algorithm. \\

Our algorithm follows the basic prototype of the stochastic gradient descent method in distributed settings. However, unlike the traditional distributed stochastic gradient descent, the server in our algorithm multiplies the stochastic gradients received from the agents by a {\em stochastic pre-conditioner} matrix. These pre-conditioned stochastic gradients are then used by the server to update the current estimate. In literature, this technique is commonly referred as {\em pre-conditioning} \cite{nocedal2006numerical}. It should be noted that unlike the conventional pre-conditioning techniques \cite{fessler2008image}, in our case, the stochastic pre-conditioning matrix itself is updated over the iterations with help from the agents. Hence, the name {\em iterative pre-conditioning}. \\

In order to present the algorithm, we introduce some notation.
The {\em individual data points} of the agents are represented by an {\em input row vector} $a$ of dimensions $1 \times d$ and a {\em scalar output} $b$. Thus, $a \in \R^{1 \times d}$ and $b \in \R$. For each data point $(a,b)$, we define {\em individual cost function} $f: \R^d \to \R$ such that for a given $x \in \R^d$,
\begin{align}
     f(x;a,b) = \frac{1}{2} \left(a ~ x - b\right)^2, \label{eqn:dist_opt}
\end{align}
and the gradient of the individual cost function $f$ as
\begin{align}
    g(x;a,b) = \nabla_x f(x;a,b) = a^T \, \left(a \, x - b\right). \label{eqn:g}
\end{align}
Here, $(\cdot)^T$ denotes the transpose. \\

In each iteration $t \in \{0, \, 1, \ldots\}$, the server maintains an estimate $x(t)$ of a minimum point~(\ref{eqn:opt_1}), and a stochastic pre-conditioner matrix $K(t) \in \R^{d \times d}$. The initial estimate $x(0)$ and the pre-conditioner matrix $K(0)$ are chosen arbitrarily from $\R^d$ and  $\R^{d \times d}$, respectively. For each iteration $t=0,1,\ldots$, the algorithm steps are presented below.

\subsection{Steps in Each Iteration $t$}
\label{sec:algo_steps}

Before initiating the iterations, the server chooses a positive scalar real-valued parameter $\beta$ and broadcast it to all the agents. We number the agents in order from $1$ to $m$. 
In each iteration $t$, the proposed algorithm comprises of four steps described below. These steps are executed collaboratively by the server and the agents, without requiring any agent to share its local data points. For each iteration $t$, the server also chooses two positive scalar real-valued parameters $\alpha$ and $\delta$.

\begin{itemize}
\setlength\itemsep{0.5em}
    \item {\em Step 1:} The server sends the estimate $x(t)$ and the pre-conditioner matrix $K(t)$ to each agent $i \in \{1,\ldots,m\}$.
    \item {\em Step 2:} 
    Each agent $i \in \{1,\ldots,m\}$ chooses a data point $(a^{i_t},b^{i_t})$ {\em uniformly} at random from its local data points $(A^i, B^i)$. Note that, $a^{i_t}$ and $b^{i_t}$ are respectively a row in the input matrix $A^i$ and the output vector $B^i$ of agent $i$. Each data point is independently and identically distributed (i.i.d.). Based on the selected data point $(a^{i_t},b^{i_t})$, each agent $i$ then computes a stochastic gradient, denoted by $g^{i_t}(t)$, which is defined as
    \begin{align}
        g^{i_t}(t) = g(x(t);a^{i_t},b^{i_t}). \label{eqn:g_i}
    \end{align}
    
    In the same step, each agent $i \in \{1,\ldots,m\}$ computes a set of vectors $\left\{R^{i_t}_j(t): ~ j = 1, \ldots, \, d \right\}$ defined as 
    \begin{align}
        R^{i_t}_j(t) = R_j(k_j(t);a^{i_t},b^{i_t}), \label{eqn:Rij}
    \end{align}
    where the function $R_j: \R^d \to \R^d$ is defined below.
    Let $I$ denote the $(d \times d)$-dimensional identity matrix. Let $e_j$ and $k_j(t)$ denote the $j$-th columns of matrices $I$ and $K(t)$, respectively. Then,
    for each column index $j \in \{1,\ldots,d\}$ of $K(t)$ and each individual data point $(a,b)$, we define
    \begin{align}
        R_j(k_j;a,b) = \left(a^T a+\beta I \right) k_j - e_j. \label{eqn:R}
    \end{align}
    
    % Each agent $l \in \{1,\ldots,m\}$ generates a realization or sample $w_l^t$ of the {\em uniform random variable} $w_l$ from the set of integers $\{1,\ldots,n\}$. Recall that, the number of local data point held by an agent is $n$. The
    % realizations $\{w_l^t, l=1,\ldots,m\}$ are independent and identically distributed (i.i.d.). Here we note that, the $w_l^t$-th data point of agent $i$ is the $(n(l-1)+w_l^t)$-th data point in the collective dataset $(A,B)$. Each agent $l$ then computes the gradient $g^{n(l-1)+w_l^t}(t)$ of the individual cost function $f^{n(l-1)+w_l^t}(x)$ based only on the randomly chosen local data point $(a^{n(l-1)+w_l^t},b^{n(l-1)+w_l^t})$.  Specifically, for each individual data point index $i\in \{1,\ldots,N\}$, 
    % \begin{align}
    %     g^i(t) = \nabla f^i(x(t)) = \left( a^i \right)^T \, \left(a^i \, x(t) - b^i\right). \label{eqn:g_i}
    % \end{align}
    % Here, $(\cdot)^T$ denotes the transpose.
    
    % In the same step, each agent $i \in \{1,\ldots,m\}$ computes a set of vectors $\left\{R^{n(l-1)+w_l^t}_j(t): ~ j = 1, \ldots, \, d \right\}$ defined as follows. Let $I$ denote the $(d \times d)$-dimensional identity matrix. Let $e_j$ and $k_j(t)$ denote the $j$-th columns of matrices $I$ and $K(t)$, respectively. Then,
    % for each column index $j$ of $K(t)$ and each individual data point index $i\in \{1,\ldots,N\}$
    % \begin{align}
    %     R^i_j(t) = \left(\left(a^{i}\right)^T a^{i}+\frac{\beta}{m} I \right) k_j(t) - \frac{1}{m} e_j. \label{eqn:Rij}
    % \end{align}
    \item {\em Step 3:} 
    Each agent $i \in \{1,\ldots,m\}$ sends the stochastic gradient $g^{i_t}(t)$ and the set of stochastic vectors
    $\left\{R^{i_t}_j(t), ~ j = 1, \ldots, \, d \right\}$ to the server.
    % Each agent $l$ sends the gradient $g^{n(l-1)+w_l^t}(t)$ and the set of vectors

    % $\left\{R^{n(l-1)+w_l^t}_j(t), ~ j = 1, \ldots, \, d \right\}$ to the server.
    \item {\em Step 4:} 
    The server draws an i.i.d. sample $\zeta_t$ {\em uniformly} at random from the set of agents $\{1,\ldots,m\}$ and updates the matrix $K(t)$ to $K(t+1)$ such that, for each $j \in \{1,\ldots,d\}$,
    \begin{align}
        k_j(t + 1) = k_j(t) - \alpha R^{\zeta_{t_t}}_j(t). \label{eqn:kcol_update_dist}
    \end{align}
    Recall that $\alpha$ is a non-negative real value. Finally, the server updates the estimate $x(t)$ to $x(t + 1)$ such that
    \begin{align}
        x(t+1) = x(t) - \delta K(t + 1) g^{\zeta_{t_t}}(t). \label{eqn:x_update_dist}
    \end{align}
    Parameter $\delta$ is a non-negative real value, commonly referred as the {\em stepsize}.
    % The server generates an i.i.d realization $k^t$ of the uniform random variable $k$ from the set of integers $\{1,\ldots,m\}$. Based on the generated number $k^t$, the server picks the corresponding agent $k^t$ and use the information received from agent $k^t$ to update the matrix $K(t)$ and the estimate $x(t)$ as follows.
    % The server updates the matrix $K(t)$ to $K(t+1)$ such that, for each $j \in \{1,\ldots,d\}$,
    % \begin{align}
    %     k_j(t + 1) = k_j(t) - \alpha R^{n(k^t-1)+w_{k^t}}_j(t). \label{eqn:kcol_update_dist}
    % \end{align}
    % Recall that $\alpha$ is a non-negative real value. Finally, the server updates the estimate $x(t)$ to $x(t + 1)$ such that
    % \begin{align}
    %     x(t+1) = x(t) - \delta K(t + 1) g^{n(k^t-1)+w_{k^t}}(t). \label{eqn:x_update_dist}
    % \end{align}
    % Parameter $\delta$ is a non-negative real value, commonly referred as the {\em step-size}.
\end{itemize}

These steps of our algorithm are summarized in Algorithm~\ref{algo_1}. 

\begin{algorithm}
  \caption{{\em The IPSG method.}}\label{algo_1}
  \begin{algorithmic}[1]
    \State The server initializes $x(0) \in \R^d$, $K(0) \in \R^{d \times d}$, $\beta > 0$ and chooses $\{\alpha > 0, \delta > 0: ~ t=0,1,\ldots\}$.
    \State {\bf Steps in each iteration} $t \in \{0, \, 1, \, 2, \ldots \}$:
    % \For{\text{each iteration $t=0, \, 1, \, 2, \ldots $}}
      \State The server sends $x(t)$ and $K(t)$ to all the agents.
      \State Each agent $i  \in \{1, \ldots, \, m\} $
      {\em uniformly} selects an i.i.d. data point $(a^{i_t},b^{i_t})$ from its local data points $(A^i,B^i)$.
      \State Each agent $i  \in \{1, \ldots, \, m\} $ sends to the server a stochastic gradient $g^{i_t}(t)$, defined in~\eqref{eqn:g_i},
      and $d$ stochastic vectors $R^{i_t}_1(t), \ldots, \, ~R^{i_t}_d(t)$, defined in~\eqref{eqn:Rij}.
      \State The server {\em uniformly} draws an i.i.d. sample $\zeta_t$ from the set of agents $\{1,\ldots,m\}$.
      \State The server updates $K(t)$ to $K(t + 1)$ as defined by~(\ref{eqn:kcol_update_dist}).
      \State The server updates the estimate $x(t)$ to $x(t + 1)$ as defined by~(\ref{eqn:x_update_dist}).
    % \EndFor
  \end{algorithmic}
\end{algorithm}

\section{CONVERGENCE ANALYSIS}
\label{sec:conv}

In this section, we present the formal convergence guarantees of Algorithm~\ref{algo_1}. We begin by introducing some notation and our main assumptions.

\subsection{Notation and Assumptions}
\label{sub:assump}

\begin{itemize}
    \item The {\em collective input matrix}, denoted by $A$, is defined to be
\begin{align}
    A = \begin{bmatrix} (A^1)^T, \ldots, \, (A^m)^T \end{bmatrix}^T, \label{eqn:data_matrix}
\end{align}
and the {\em collective output vector}, denoted by $B$, is defined to be
\begin{align}
    B = \begin{bmatrix} (B^1)^T, \ldots, \, (B^m)^T \end{bmatrix}^T. \label{eqn:b_vector}
\end{align}
Define $N = mn$.
Note that, $A \in \R^{N \times d}$ and $B \in \R^{N}$.
    \item Let $\nabla F(t)$ denote the true gradient of the objective cost $\frac{1}{m}\sum_{i=1}^m F^i$ evaluated at the current estimate $x(t)$ of the solution in~\eqref{eqn:opt_1}. 
    \item Since the matrix $\A$ is symmetric positive semi-definite, its eigenvalues are real non-negative. We let $s_1 \geq \ldots \geq s_d \geq 0$ denote the eigenvalues of $\A$. For $\beta > 0$, we define
\begin{align}
    K_{\beta} = \left(\frac{1}{N}\A + \beta I\right)^{-1}. \label{eqn:kbeta}
\end{align}
Since $\beta > 0$ and $\A$ is positive semi-definite, $K_{\beta}$ is well-defined.
    \item For each iteration $t\geq 0$ and agent $i \in \{1,\ldots,m\}$, let $\condexp{\cdot}$ denote the conditional expectation of a function the random variable $i_t$ given the current estimate $x(t)$ and the current pre-conditioner $K(t)$. \item For each iteration $t\geq 0$, let $\serverexp{\cdot}$ denote the conditional expectation of a function the random variable $\zeta_t$ given the current estimate $x(t)$ and the current pre-conditioner $K(t)$.
    \item For each iteration $t\geq 0$, let 
    \begin{align}
        I_t = \{i_t, ~ i=1,\ldots,m\} \cup \{\zeta_t\}. \label{eqn:It}
    \end{align}
    and
    \begin{align}
        \agentexp{\cdot} = \E_{1_t, \ldots m_t,\zeta_t}(\cdot). \label{eqn:agentexp_def}
    \end{align}
    \item Let $\totexp{\cdot}$ denote the total expectation of a function of the collection of the random variables $\{I_0,\ldots,I_t\}$ given the initial estimate $x(0)$ and initial pre-conditioner matrix $K(0)$. Specifically,
    \begin{align}
        \totexp{\cdot} = \E_{I_0,\ldots,I_t}(\cdot), ~ t\geq 0. \label{eqn:def_totexp}
    \end{align}
    % \item Let $\expect{\cdot}$ denote the expectation of a function of uniform random variable in $\{1,\ldots,N\}$.
    \item For each iteration $t\geq 0$, define the conditional variance of the stochastic gradient $g^{i_t}(t)$, which is a function of the random variable $i_t$, given the current estimate $x(t)$ and the current pre-conditioner $K(t)$ as
    \begin{align}
    \var{g^{i_t}(t)} & = \condexp{\norm{g^{i_t}(t) - \condexp{g^{i_t}(t)}}^2} \nonumber \\
    & = \condexp{\norm{g^{i_t}(t)}^2} - \norm{\condexp{g^{i_t}(t)}}^2. \label{eqn:variance} 
    \end{align}
\end{itemize}

We make the following assumption on the rank of the matrix $\A$.

\begin{assumption}
\label{asp:full_rank}
Assume that the matrix $\A$ is full rank. 
\end{assumption}
% \noindent \textbf{Assumption~1}: 

Note that Assumption~\ref{asp:full_rank} holds true if and only if the matrix $\A$ is positive definite with $s_d > 0$. 
As the Hessian of the aggregate cost function $\sum_{i=1}^m F^i(x)$ is equal to $\A$ for all $x$ (see~\eqref{eqn:dist}, under Assumption~\ref{asp:full_rank}, the aggregate cost function has a unique minimum point. Equivalently, the solution of the distributed least squares problem defined by~\eqref{eqn:opt_1} is unique. \\

We also assume, as formally stated in Assumption~\ref{asp:bnd_var} below, that the variance of the stochastic gradient for each agent is bounded. This is a standard assumption for the analysis of stochastic algorithms~\cite{bottou2018optimization}. 

\begin{assumption}
\label{asp:bnd_var}
    Assume that there exist two non-negative real scalar values $V_1$ and $V_2$ such that, for each iteration $t=0,1,\ldots$ and each agent $i \in \{1,\ldots,m\}$, 
    \begin{align*}
        \var{g^{i_t}(t)} \leq V_1 + V_2\norm{\nabla F(t)}^2.
    \end{align*}
\end{assumption}

Next, we present our key result on the convergence of Algorithm~\ref{algo_1}.
% is presented below in Theorem~\ref{thm:z_conv}. 

\subsection{Main Theorem}
\label{sub:theorem}

To formally present our key result in Theorem~\ref{thm:z_conv} below, we introduce some additional notation.
% We require a few more notation for presenting our key result.
\begin{itemize}
    \item Let $a^i$ and $b^i$, respectively, denote the $i$-th row of the input matrix $A$ and output vector $B$.
    
    \item Let $\Lambda_{i}$ and $\lambda_{i}$, respectively, denote the largest and the smallest eigenvalue of the positive semi-definite matrix $\Ai$.
    
    \item Recall the definition of $K_{\beta}$ in~\eqref{eqn:kbeta}. For positive real values $\alpha$ and $\beta$, we define
\begin{align}
    C_2 & =  \alpha ~ \frac{1}{N} \sum_{i=1}^N \norm{\Ai-\frac{1}{N}\A} \norm{K_{\beta}}. \label{eqn:c2} 
\end{align}
    \item We introduce the following parameter $\rho$ which determines the rate of convergence of the pre-conditioner matrix $K(t)$. For positive real values $\alpha$ and $\beta$, we define
\begin{align}
    \rho & = \frac{1}{N} \sum_{i=1}^N \norm{I - \alpha \left(\Ai+\beta I \right)}. \label{eqn:rho}
\end{align}
    
    \item For the positive valued parameter $\beta$, let
\begin{align}
    L & = \beta + \max_{i=1,\ldots,N} \Lambda_i, \label{eqn:L} \text{ and } \\
    \sigma^2 & = \max_{j=1,\ldots,d} \frac{1}{N} \sum_{i=1}^N \norm{\left(\Ai+\beta I \right) K_{\beta}~ e_j  - e_j}^2. \label{eqn:sigma}
\end{align}
    \item Recall from Section~\ref{sub:assump} that $s_d > 0$ is the smallest eigen value of the positive definite matrix $\A$. For the positive valued parameters $\alpha$ and $\beta$, let
\begin{align}
    C_3 & = \frac{\alpha N \sigma^2}{s_d \left(1-\alpha L\right)}. \label{eqn:c3}
\end{align}
    \item We define the following parameter $\overline{\alpha}$, determining a sufficient range of the parameter $\alpha$ in Algorithm~\ref{algo_1}:
\begin{align}
    \Bar{\alpha} & = \min \left\{\frac{N}{s_d}, \frac{1}{L}, \frac{2}{(s_1/N)+\beta}\right\}. \label{eqn:alpha_bar}
\end{align}
% In order to prove convergence of Algorithm~\ref{algo_1}, for each iteration $t$, we introduce a sufficient range $\overline{\delta}(t)$ for the stepsize $\delta$. This parameter $\overline{\delta}(t)$ has been defined latter in~\eqref{eqn:delta_bar} (see Appendix~\ref{sub:notations}).
    \item For each iteration $t\geq 0$, and the positive valued parameters $\alpha$, $\beta$ and $\delta$, we define the following parameter $R_1(t)$ which determines the rate of convergence of Algorithm~\ref{algo_1}:
\begin{align}
    R_1(t) & = 1 + \delta^2 C_8(t) +\alpha \delta C_5(t) - \delta C_6(t). \label{eqn:r1}
\end{align}
The parameters $C_5(t)$, $C_6(t)$ and $C_8(t)$ are defined latter in Appendix~\ref{sub:notations}. Note that $C_5(t)$, $C_6(t)$ and $C_8(t)$ depend on $\alpha$, and are independent of the step-size $\delta$.
    \item For each iteration $t\geq 0$, we define the estimation error
\begin{align}
    z(t) = x(t) - x^*, ~ t\geq 0. \label{eqn:err}
\end{align}
\end{itemize}

\begin{theorem} \label{thm:z_conv}
Consider Algorithm~\ref{algo_1} with parameters $\beta > 0$, $\alpha < \overline{\alpha}$ and $\delta > 0$. If Assumptions~\ref{asp:full_rank} and~\ref{asp:bnd_var} are satisfied, then there exist two non-negative real scalar values $E_1 \geq \sqrt{V_1N}$ and $E_2 \geq \sqrt{V_2N}$ such that the following statements hold true.
\begin{enumerate}[(i)]
    \item If the stepsize $\delta$ is sufficiently small, then there exists a non-negative integer $T < \infty$ such that for any iteration $t \geq T$, $R_1(t)$ is positive and less than $1$.
    
    \item For an arbitrary time step $t \geq 0$, given the estimate $x(t)$ and the matrix $K(t)$, 
    % for all $t \geq 0$,
    \begin{align}
        \totexp{\norm{z(t+1)}^2} \leq R_1(t) \norm{z(t)}^2 + R_2(t). \label{eqn:iter_err}
    \end{align}
    
    \item Given arbitrary choices of the initial estimate $x(0) \in \R^d$ and the pre-conditioner matrix $K(0) \in \R^{d \times d}$,
    \begin{align}
        \lim_{t \to \infty} \totexp{\norm{z(t+1}^2}
        \leq \delta^2 V_1 N \left(d C_3 + \norm{K_{\beta}}^2 + \frac{2 C_2 \norm{K_{\beta}}}{1-\rho}  \right) + 2 \alpha^2 \left(C_1 E_1 \norm{K_{\beta}} \right)^2. \label{eqn:final_err}
    \end{align}
\end{enumerate}
\end{theorem}
~\\

The proof of Theorem~\ref{thm:z_conv} is deferred to Appendix~\ref{prf:z_conv}.\\

The implications of Theorem~\ref{thm:z_conv} are as follows.
\begin{itemize}
    \item According to Part (i) and (ii) of Theorem~\ref{thm:z_conv}, for small enough values of the parameters $\alpha$ and stepsize $\delta$, as $R_1(t) \in (0,1)$ after a finite number of iterations, Algorithm~\ref{algo_1} converges {\em linearly} in expectation to a neighborhood of the minimum point $x^*$ of the distributed least-squares problem~\eqref{eqn:opt_1}.
    \item According to Part (iii) of Theorem~\ref{thm:z_conv}, the neighborhood of $x^*$, to which the estimates of Algorithm~\ref{algo_1} converges in expectation, is $\bigO(V_1)$. In other words, the sequence of expected ``distance'' between the minima $x^*$ of~\eqref{eqn:opt_1} and the final estimated value of Algorithm~\ref{algo_1} is proportional to the variance of the stochastic gradients at the minimum point.
    % \item A choice of vanishing stepsizes $\delta$ such that $\delta \to 0$ as $t \to \infty$ ensures zero error in the final estimated value, as evident from~\eqref{eqn:final_err}. Although, the convergence will no longer be linear in that case.
    \item For the special case, when the same minimum point $x^*$ minimizes each of the individual cost function $f^i$ (although $x^*$ need not be the unique minimizer for any of them), we have $V_1 = 0$ and hence $\sigma^2 = 0$ (see~\eqref{eqn:sigma}). In that case,~\eqref{eqn:final_err} implies that the sequence of estimates $\{x(t), ~ t=0,1,\ldots\}$ in Algorithm~\ref{algo_1} converges exactly to the minimum point.
\end{itemize}

% \textcolor{red}{ \textbf{*********************************} } \\

% As $\rho \in [0,1)$,~\eqref{eqn:conv_1} of Theorem~\ref{thm:z_conv} implies that
% \begin{align}
%     \underset{t \rightarrow \infty}{\lim} \frac{\norm{\Et \left[z(t+1)\right]}}{\norm{z(t)}} \leq \mu < 1. \label{eqn:conv_rate}
% \end{align}
% The above inequality signifies that the sequence of expected value of the estimation error $\{\Et \left[z(t)\right]\}_{t \geq 0}$ converges linearly to $0_d$ with {\em rate of convergence} equal to $\mu$. Equivalently, the sequence of expected value of the estimates $\{\Et \left[x(t)\right]\}_{t \geq 0}$ linearly converges to a minimum point $x^*$ of the expected cost function $\E_{i\sim \D} [F^i(x)]$ in~\eqref{eqn:opt_1}.

\section{EXPERIMENTAL RESULTS}
\label{sec:exp}

In this section, we present our experimental results validating the convergence of our proposed algorithm on real-world problems and its comparison with related methods.  

\subsection{Setup}

We conduct experiments for different collective data points $(A,B)$. These collective input matrices $A$ have been chosen such that the set of their condition number covers a wide range of values. Four of these datasets are from the the benchmark datasets available in SuiteSparse Matrix Collection\footnote{https://sparse.tamu.edu}. Particularly these four datasets are \textit{``ash608''}, \textit{``abtaha1''}, \textit{``gre\_343''}, and \textit{``illc1850''}. The fifth dataset, \textit{``cleveland''}, is from the UCI Machine Learning Repository~\cite{Dua:2019}. The sixth and final dataset is the \textit{``MNIST''}\footnote{https://www.kaggle.com/oddrationale/mnist-in-csv} dataset. \\

In the case of the first four aforementioned datasets, the problem is set up as follows. Consider a particular dataset \textit{``ash608''}. Here, the matrix $A$ has $608$ rows and $d=188$ columns. The collective output vector $B$ is such that $B=Ax^*$ where $x^*$ is a $188$ dimensional vector, all of whose entries are unity. To simulate the distributed server-agent architecture, the data points represented by the rows of the matrix $A$ and the corresponding observations represented by the elements of the vector $B$ are divided amongst $m=8$ agents numbered from $1$ to $8$. Since the matrix $A$ for this particular dataset has $608$ rows and $188$ columns, each of the eight agents $1, \ldots, \, 8$ has a data matrix $A^i$ of dimension $76 \times 188$ and a observation vector $B^i$ of dimension $76$. 
The data points for the other three datasets, \textit{``abtaha1''}, \textit{``gre\_343''}, and \textit{``illc1850''}, are similarly distributed among $m=4$, $m=7$ and $m=10$ agents, respectively.\\

For the fifth dataset, $212$ arbitrary instances from the \textit{``cleveland''} dataset have been selected. This dataset has $13$ measured or derived numeric attributes, each corresponding to a column in the matrix $A$, and a target class (whether the patient has heart disease or not), which corresponds to the output vector $B$. Since the attributes in the matrix $A$ has different units, each column in $A$ is then shifted by the mean value of the entries in the corresponding column and then divided by the standard deviation of the entries in that column. Finally, a $212$-dimensional column vector of unity is appended to this pre-processed matrix. This is our final input matrix $A$ of dimension $(212 \times 14)$. The collective data points $(A,B)$ are then distributed among $m=4$ agents, in the manner described earlier.\\

From the training examples of the \textit{``MNIST''} dataset, we consider those instances labeled as either the digit one or the digit five and select $1500$ such rows arbitrarily. For each instance, we calculate two quantities. One is the average intensity of an image, defined as the sum of the pixels divided by the number of pixels in each image. The other one is the average symmetry of an image, which is defined as the negation of the average absolute deviation of the pixel values between an image and its flipped version. Let the average intensity of all $1500$ images be denoted by the $1500$-dimensional column vector $a_1$. Similarly, let the average symmetry of those $1500$ images be denoted by the $1500$-dimensional column vector $a_2$. Then, our input matrix before pre-processing is $\begin{bmatrix} a_1 & a_2 & a_1.^2 & a_1.*a_2 & a_2.^2 \end{bmatrix}$. Here, $(.*)$ represents element-wise multiplication and $(.^2)$ represents element-wise squares. This raw input matrix is then pre-processed as described earlier for the \textit{``cleveland''} dataset. Finally, a $1500$-dimensional column vector of unity is appended to this pre-processed matrix. This is our final input matrix $A$ of dimension $(1500 \times 6)$. The collective data points $(A,B)$ are then distributed among $m=10$ agents, in the manner already described for the other datasets.\\

As the matrix $\A$ is positive definite in each of these cases, the optimization problem~(\ref{eqn:opt_1}) has a unique solution $x^*$ for all of these datasets.

\begin{table*}[htb!]
\caption{\it The parameters used in different algorithms.}
\begin{center}
\begin{tabular}{|p{1.4cm}||p{2cm}|p{1.7cm}|p{2.4cm}|p{2.5cm}|p{2.5cm}|}
\hline
Dataset & IPSG & SGD~\cite{bottou2018optimization} & AdaGrad~\cite{duchi2011adaptive} & AMSGrad~\cite{reddi2019convergence} & Adam~\cite{kingma2014adam} \\
\hline
\hline
cleveland & $\alpha = 0.0031$, $\delta = 0.5$, $\beta = 30$ & $\alpha = 0.0031$ & $\alpha = 1$, $\epsilon = 10^{-7}$ & $\alpha = 0.05$, $\beta_1 = 0.9$, $\beta_2 = 0.999$, $\epsilon = 10^{-7}$ & $\alpha = 0.05$, $\beta_1 = 0.9$, $\beta_2 = 0.999$, $\epsilon = 10^{-7}$ \\
\hline
ash608 & $\alpha = 0.1163$, $\delta = 1$, $\beta = 1$ & $\alpha = 0.1163$ & $\alpha = 1$, $\epsilon = 10^{-7}$ & $\alpha_t = \frac{0.5}{\sqrt{t}}$, $\beta_1 = 0.9$, $\beta_2 = 0.99$, $\epsilon = 10^{-7}$ & $\alpha_t = \frac{0.1}{\sqrt{t}}$, $\beta_1 = 0.9$, $\beta_2 = 0.999$, $\epsilon = 10^{-7}$ \\
\hline
abtaha1 & $\alpha = 0.0052$, $\delta = 2$, $\beta = 5$ & $\alpha = 0.0052$ & $\alpha = 1$, $\epsilon = 10^{-7}$ & $\alpha_t = \frac{1}{\sqrt{t}}$, $\beta_1 = 0.9$, $\beta_2 = 0.99$, $\epsilon = 10^{-7}$ & $\alpha_t = \frac{0.5}{\sqrt{t}}$, $\beta_1 = 0.9$, $\beta_2 = 0.999$, $\epsilon = 10^{-7}$ \\
\hline
MNIST & $\alpha = 0.0003$, $\delta = 0.1$, $\beta = 1$ & $\alpha = 0.0003$ & $\alpha = 1$, $\epsilon = 10^{-7}$ & $\alpha = 1$, $\beta_1 = 0.9$, $\beta_2 = 0.999$, $\epsilon = 10^{-7}$ & $\alpha = 0.1$, $\beta_1 = 0.9$, $\beta_2 = 0.999$, $\epsilon = 10^{-7}$ \\
\hline
gre\_343 & $\alpha = 1.2$, $\delta = 2.5$, $\beta = 0.5$ & $\alpha = 1.96$ & $\alpha = 1$, $\epsilon = 10^{-7}$ & $\alpha_t = \frac{0.1}{\sqrt{t}}$, $\beta_1 = 0.9$, $\beta_2 = 0.999$, $\epsilon = 10^{-7}$ & $\alpha_t = \frac{0.2}{\sqrt{t}}$, $\beta_1 = 0.9$, $\beta_2 = 0.999$, $\epsilon = 10^{-7}$ \\
\hline
illc1850 & $\alpha = 0.4436$, $\delta = 2$, $\beta = 1$ & $\alpha = 0.4436$ & $\alpha = 1$, $\epsilon = 10^{-7}$ & $\alpha_t = \frac{0.5}{\sqrt{t}}$, $\beta_1 = 0.9$, $\beta_2 = 0.99$, $\epsilon = 10^{-7}$ & $\alpha_t = \frac{0.5}{\sqrt{t}}$, $\beta_1 = 0.9$, $\beta_2 = 0.999$, $\epsilon = 10^{-7}$ \\
\hline
\end{tabular}
\end{center}
\label{tab:parameters}
\end{table*}

\subsection{Comparisons with Existing Algorithms}

We compare the performance of our proposed algorithm on the aforementioned datasets, with the other stochastic algorithms mentioned in Section~\ref{sec:intro}. Specifically, these algorithms are stochastic gradient descent (SGD)~\cite{bottou2018optimization}, adaptive gradient descent (AdaGrad)~\cite{duchi2011adaptive}, adaptive momentum estimation (Adam)~\cite{kingma2014adam}, and AMSGrad~\cite{reddi2019convergence} in the distributed network architecture of Fig.~\ref{fig:sys}.

\subsubsection{Algorithm Parameter Selection}

These algorithms are implemented with different combinations of the respective algorithm parameters on the individual datasets. The parameter combinations are described below.
\begin{enumerate}[(i)]
    \item {\bf IPSG}: The optimal (smallest) convergence rate of the deterministic version of the proposed IPSG method (Algorithm~\ref{algo_1}) is obtained when $\alpha = \frac{2}{s_1 + s_d}$~\cite{chakrabarti2020iterative}. For each of the six datasets, we find that the IPSG method converges fastest when the parameter $\alpha$ in the proposed IPSG algorithm is set similarly as $\alpha = \frac{2}{s_1 + s_d}$. Recall that, $s_1$ and $s_d$ are respectively the largest and the smallest eigenvalue of the positive definite matrix $\A$. The stepsize parameter $\delta$ of the IPSG algorithm is chosen from the set $\{0.1,0.5,1,2,2.5\}$. The parameter $\beta$ is chosen from the set $\{0.1,0.5,1,5,10,30,50\}$.
    \item {\bf SGD}: The SGD algorithm has only one parameter: the stepsize, denoted as $\alpha$~\cite{bottou2018optimization}. We choose $\alpha$ for SGD following the same reasoning as done for the IPSG method above. Specifically, the deterministic version of the SGD method is the gradient-descent method, which has the optimal rate of convergence when $\alpha = \frac{2}{s_1 + s_d}$. We find that the SGD method converges fastest when the stepsize parameter is similarly set as $\alpha = \frac{2}{s_1 + s_d}$. It must be noted that $s_1$ and $s_d$ depends on the collective data matrix $A$, and hence their values may not be known to the server. When the actual values or estimates of $s_1$ and $s_d$ are not known, the parameter $\alpha$ in both the IPSG and the SGD algorithm can be experimentally set by trying several values of different orders, as done for the parameters $\delta$ and $\beta$ in the IPSG method above.
    \item {\bf AdaGrad}: The stepsize parameter $\alpha$ of the AdaGrad algorithm~\cite{duchi2011adaptive} is selected from the set $\{0.1,1,\frac{1}{t}\}$. The parameter $\epsilon$ is set at its usual value of $10^{-7}$.
     \item {\bf Adam} and {\bf AMSGrad}: The stepsize parameter $\alpha$ of these two algorithms~\cite{kingma2014adam, reddi2019convergence} is selected from the set $\{c,\frac{c}{\sqrt{t}}\}$ where $c$ is from the set $\{1,0.5,0.1,0.2,0.05,0.01\}$. The parameter $\beta_1$ is set at its usual value of $0.9$. The parameter $\beta_2$ is selected from their usual values of $\{0.99,0.999\}$. The parameter $\epsilon$ is set at $10^{-7}$.
\end{enumerate}
The best parameter combination from above, for which the respective algorithms converge in a fewer number of iterations, is reported for each dataset in Table~\ref{tab:parameters}. \\

\textbf{Initialization:}
The initial estimate $x(0)$ for all of these algorithms is chosen as the $d$-dimensional zero vector for each dataset except the \textit{``cleveland''} dataset. For \textit{``cleveland''} dataset, $x(0)$ is chosen as the $d$-dimensional vector whose each entry is $10$. The initial pre-conditioner matrix $K(0)$ for the IPSG algorithm is the zero matrix of dimension $(d \times d)$.

\begin{figure*}[htb!]
\centering
\begin{subfigure}{.5\textwidth}
  \begin{center}
  \includegraphics[width = \textwidth]{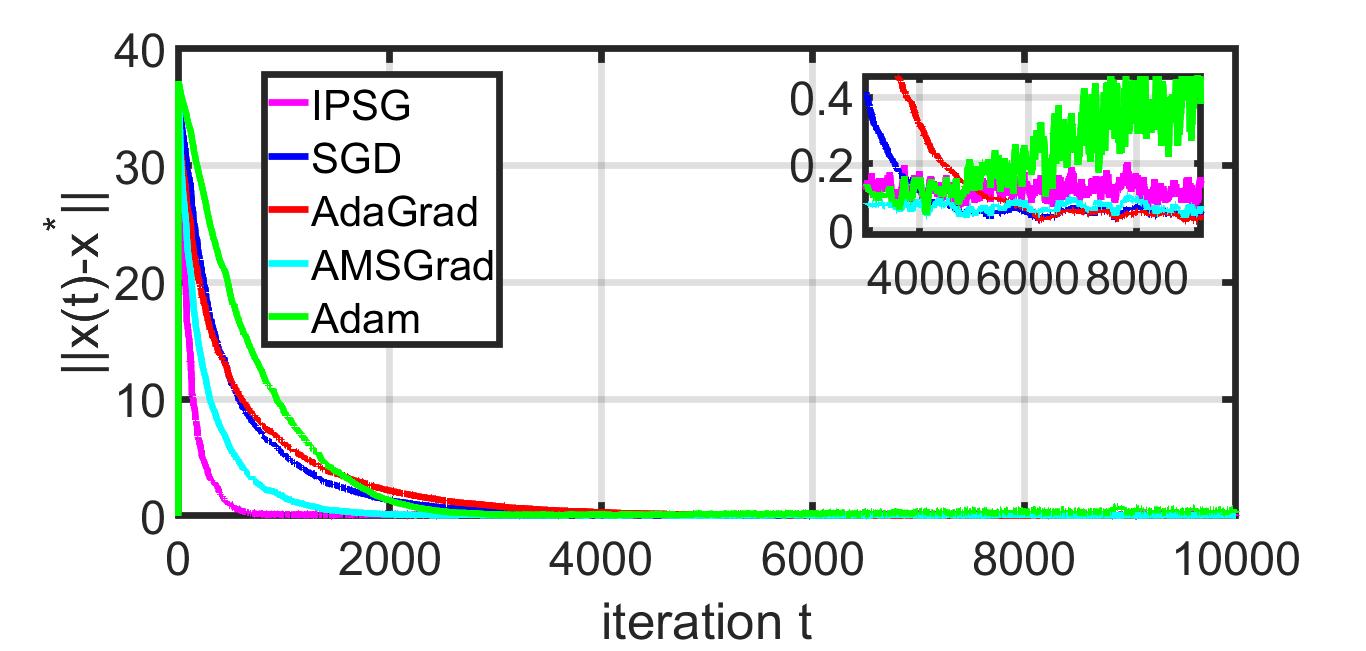}
  \caption{\textit{``cleveland''}}
  \end{center}
\end{subfigure}%
\begin{subfigure}{.5\textwidth}
  \begin{center}
  \includegraphics[width = \textwidth]{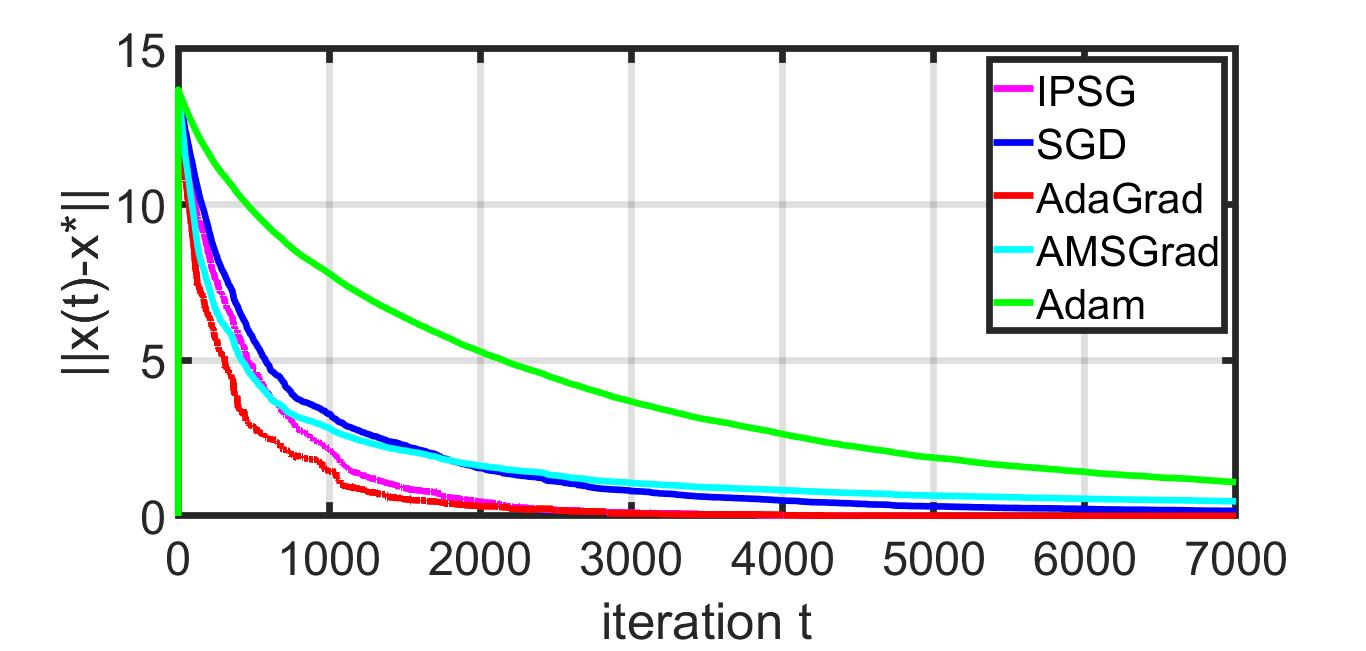}
  \caption{\textit{``ash608''}}
  \end{center}
\end{subfigure}
\begin{subfigure}{.5\textwidth}
  \begin{center}
  \includegraphics[width = \textwidth]{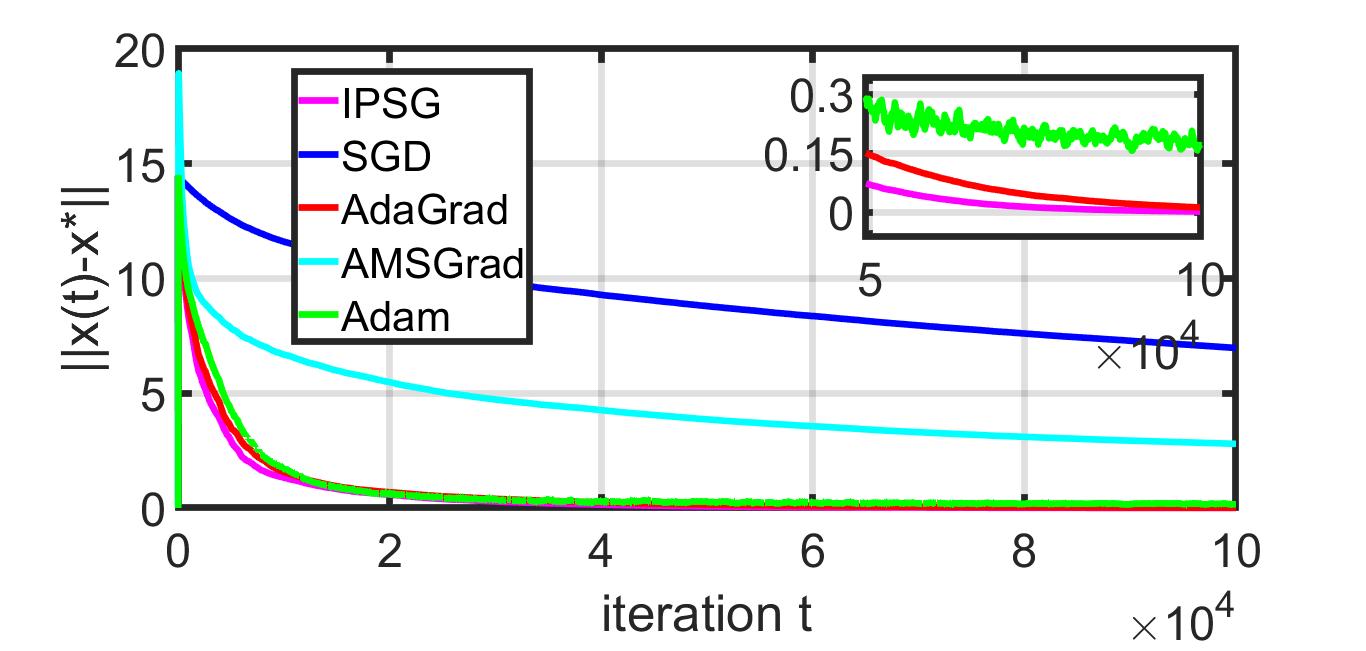}
  \caption{\textit{``abtaha1''}}
  \end{center}
\end{subfigure}%
\begin{subfigure}{.5\textwidth}
  \begin{center}
  \includegraphics[width = \textwidth]{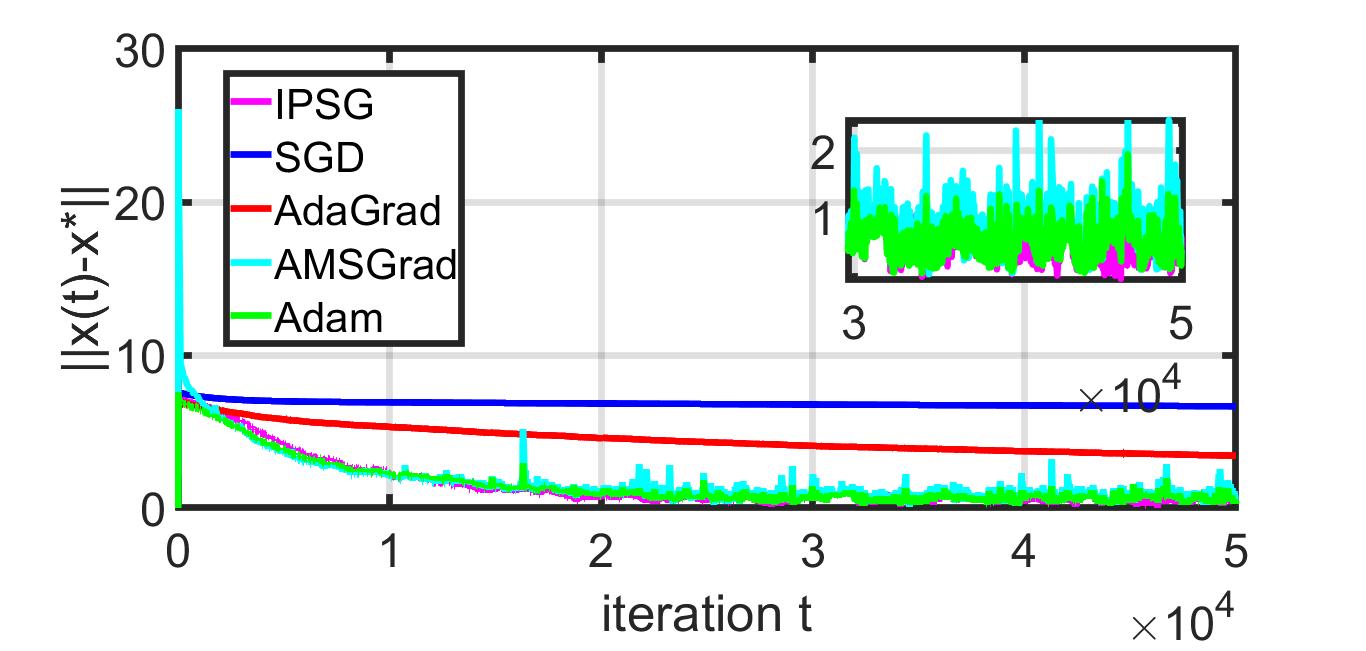}
  \caption{\textit{``MNIST''}}
  \end{center}
\end{subfigure}
\begin{subfigure}{.5\textwidth}
  \begin{center}
  \includegraphics[width = \textwidth]{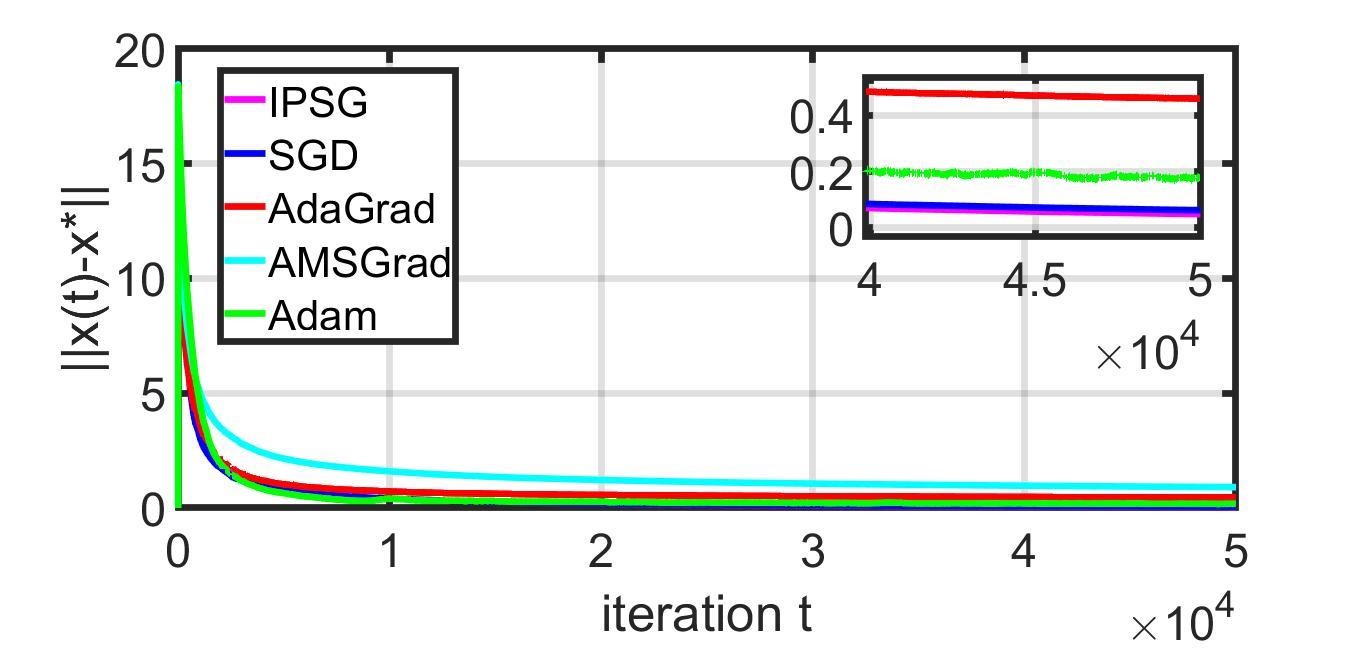}
  \caption{\textit{``gre_343''}}
  \end{center}
\end{subfigure}%
\begin{subfigure}{.5\textwidth}
  \begin{center}
  \includegraphics[width = \textwidth]{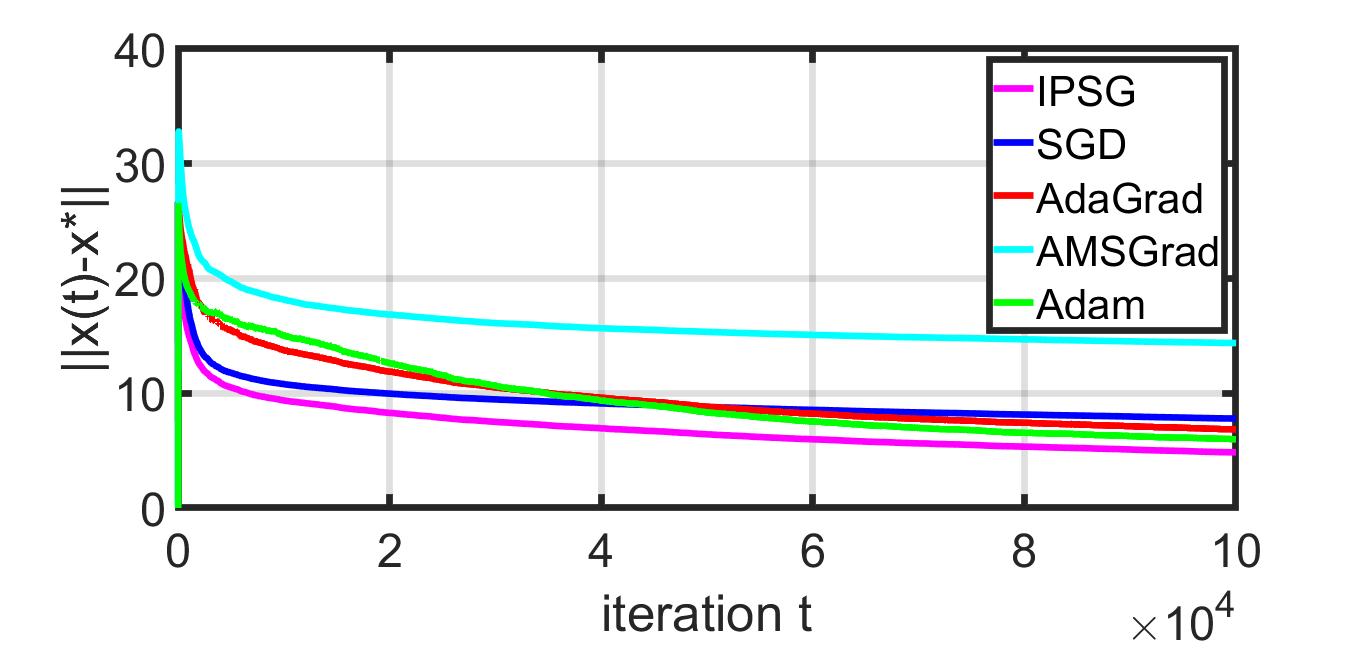}
  \caption{\textit{``illc1850''}}
  \end{center}
\end{subfigure}
\caption{\footnotesize{\it Temporal evolution of estimation error $\norm{x(t)-x^*}$, for different algorithms represented by different colors. For all the algorithms, (a): $x(0) = [10,\ldots,10]^T$; (b)-(f): $x(0) = [0,\ldots,0]^T$.
Additionally, for IPSG, $K(0) = O_{d \times d}$. The other parameters are enlisted in Table~\ref{tab:parameters}.}}
\label{fig:comp}
\end{figure*}

\begin{table*}[htb!]
\caption{\it Comparisons between the number of iterations required by different algorithms to attain the specified values for the relative estimation errors $\epsilon_{tol} = \norm{x(t)-x^*}/\norm{x(0)-x^*}$.}
\begin{center}
\begin{tabular}{|p{1.4cm}|p{1.5cm}|p{1.5cm}||p{1.5cm}|p{1.5cm}|p{1.5cm}|p{1.6cm}|p{1.5cm}|}
\hline
Dataset & $\kappa(\A)$ & $\epsilon_{tol}$ & IPSG & SGD & AdaGrad & AMSGrad & Adam \\
\hline
\hline
cleveland & $7.34$ & $1.5 \times 10^{-3}$ & $4.11 \times 10^3$ & $4.71 \times 10^3$ & $6.04 \times 10^3$ & \cellcolor{lightgray} $3.63 \times 10^3$ & $4.11 \times 10^3$ \\
\hline
ash608 & $11.38$ & $10^{-4}$ & $ \cellcolor{lightgray} 5.73 \times 10^3$ & $2.1 \times 10^4$ & $5.86 \times 10^3$ & $>4 \times 10^4$ & $>4 \times 10^4$ \\
\hline 
abtaha1 & $1.5 \times 10^2$ & $10^{-3}$ & $ \cellcolor{lightgray} 7.35 \times 10^4$ & $>10^5$ & $9.75 \times 10^4$ & $>10^5$ & $>10^5$ \\
\hline
MNIST & $2.59 \times 10^3$ & $2.6 \times 10^{-3}$ & $ \cellcolor{lightgray} 3.41 \times 10^4$ & $>5 \times 10^4$ &  $>5 \times 10^4$ & $>5 \times 10^4$ & $4.41 \times 10^4$ \\
\hline
gre\_343 & $1.25 \times 10^4$ & $4 \times 10^{-3}$ & $ \cellcolor{lightgray} 3.88 \times 10^4$ & $4.43 \times 10^5$ & $>10^5$ & $>10^5$ & $>10^5$ \\
\hline
illc1850 & $1.93 \times 10^6$ & $0.2$ & $ \cellcolor{lightgray} 8.06 \times 10^4$ & $3.31 \times 10^5$ & $2.81 \times 10^5$ & $>5 \times 10^5$ & $1.63 \times 10^5$ \\
\hline
\end{tabular}
\end{center}
\label{tab:time_comp}
\end{table*}

\subsubsection{Results}

We compare the number of iterations needed by these algorithms to reach a {\em relative estimation error} defined as 
\begin{eqnarray*}
    \epsilon_{tol} = \frac{\norm{x(t) - x^*}}{\norm{x(0) - x^*}}.
\end{eqnarray*}
Each iterative algorithm is run until its relative estimation error does not exceed $\epsilon_{tol}$ over a period of $10$ consecutive iterations, and the smallest such iteration is reported.
The algorithm parameter choices for this comparison have been described above, and their specific values are tabulated in Table~\ref{tab:parameters} for all six datasets. The comparison results are in Table~\ref{tab:time_comp}. The second column of  Table~\ref{tab:time_comp} indicates the condition number of the matrix $\A$ for each dataset. The condition numbers of these considered datasets range from a small value of $7.34$ to a large value of order $10^6$. From Table~\ref{tab:time_comp}, we see that the IPSG algorithm converges fastest among the algorithms on all the datasets but the \textit{``cleveland''} one. Note that the \textit{``cleveland''} dataset has the smallest condition number of $7.34$ among the datasets. Even for this \textit{``cleveland''} dataset, only the AMSGrad algorithm requires fewer iterations than IPSG. Moreover, from the datasets \textit{``MNIST''}, \textit{``gre\_343''}, and \textit{``illc1850''}, we observe that the differences between the proposed IPSG method and the other methods are significant when the condition number of the matrix $\A$ is of order $10^3$ or more. Thus, our claim on improvements over the prominent stochastic algorithms for solving the distributed least-squares problem~\eqref{eqn:opt_1} is corroborated by the above experimental results.

% \begin{table*}[htb!]
% \caption{\it Comparisons between the number of iterations required by different algorithms to attain the specified values for the relative estimation errors $\epsilon_{tol} = \norm{x(t)-x^*}/\norm{x(0)-x^*}$.}
% \begin{center}
% \begin{tabular}{|p{1.3cm}|p{1.5cm}|p{1.5cm}||p{1.5cm}|p{1.5cm}|p{1.5cm}|p{1.5cm}|p{1.5cm}|}
% \hline
% Dataset & $\kappa(\A)$ & $\epsilon_{tol}$ & IPSG & SGD & AdaGrad & AMSGrad & Adam \\
% \hline
% \hline
% cleveland & $7.34$ & $2\times 10^{-3}$ & $1.03 \times 10^3$ & $1.15 \times 10^3$ & $1.55 \times 10^3$ & \cellcolor{lightgray} $8.43 \times 10^2$ & $1.88 \times 10^3$ \\
% \hline
% ash608 & $11.38$ & $10^{-4}$ & $ \cellcolor{lightgray} 6.41 \times 10^2$ & $2.7 \times 10^3$ & $7.93\times 10^2$ & $3.92 \times 10^3$ & $5.99 \times 10^3$ \\
% \hline 
% abtaha1 & $1.5 \times 10^2$ & $10^{-4}$ & $ \cellcolor{lightgray} 2.64 \times 10^4$ & $>10^5$ & $3.45 \times 10^4$ & $>10^5$ & $>10^5$ \\
% \hline
% MNIST & $2.59 \times 10^3$ & $2 \times 10^{-2}$ & $ \cellcolor{lightgray} 2.82 \times 10^3$ & $>5 \times 10^4$ &  $>5 \times 10^4$ & $>5 \times 10^4$ & $6.21 \times 10^3$ \\
% \hline
% gre\_343 & $1.25 \times 10^4$ & $5 \times 10^{-4}$ & $ \cellcolor{lightgray} 4.37 \times 10^4$ & $6.28 \times 10^4$ & $>10^5$ & $>10^5$ & $>10^5$ \\
% \hline
% illc1850 & $1.93 \times 10^6$ & $10^{-1}$ & $ \cellcolor{lightgray} 6.99 \times 10^4$ & $2.52 \times 10^5$ & $2.95 \times 10^5$ & $>4 \times 10^5$ & $1.37 \times 10^5$ \\
% \hline
% \end{tabular}
% \end{center}
% \label{tab:time_comp}
% \end{table*}

\section{SUMMARY}
\label{sec:summary}

In this paper, we have considered solving the multi-agent distributed linear least-squares problem in server-based network architecture. We have proposed an iterative algorithm, namely the Iteratively Pre-Conditioned Stochastic Gradient-descent (IPSG) method, wherein the estimate of the solution is updated by the server using a single randomly sampled data point, held by the agents, at each iteration. After presenting the IPSG method in detail, we have rigorously analyzed its convergence. Subsequently, we have presented our experimental results on six real-world datasets, including the well-known MNIST dataset. The experimental results have reinforced our claim on the proposed IPSG method's superior rate of convergence compared to the state-of-the-art stochastic algorithms. We can reach a satisfactory neighborhood of the desired solution in a fewer number of iterations compared to the existing state-of-the-art stochastic algorithms.

%%%%%%%%%%%%%%%%% REFERENCES %%%%%%%%%%%%%%%%%%%%%%%% 
\newpage
\bibliographystyle{unsrt}
\bibliography{refs}

%%%%%%%%%%%%%%%%% APPENDIX %%%%%%%%%%%%%%%%%%%%%%%%%%%
\newpage
\appendix
\section{PROOFS}

\subsection{Preliminary Results}
\label{sub:prelim}

The results below are used in the proof of our main result, i.e., Theorem~\ref{thm:z_conv}.\\

\textbf{(a)}
Consider an arbitrary iteration $t\geq 0$.
Upon taking conditional expectation $\agentexp{\cdot}$ on both sides of~\eqref{eqn:Kt_iter}, given the current matrix $K(t)$ and estimate $x(t)$, we have
\begin{align}
    \agentexp{\widetilde{K}(t+1)} = \left(I - \alpha \left(\agentexp{\Azi}+\beta I \right)\right)\widetilde{K}(t) - \alpha \left(\agentexp{\Azi}-\frac{1}{N}\A \right)K_{\beta}. \label{eqn:Kt_iter_exp}
\end{align}
Upon subtracting both sides of~\eqref{eqn:Kt_iter} from that of~\eqref{eqn:Kt_iter_exp} we get
\begin{align}
    \agentexp{\widetilde{K}(t+1)} - \widetilde{K}(t+1) & = \alpha \left(\Azi - \agentexp{\Azi}\right) \left( \widetilde{K}(t) + K_{\beta} \right) \nonumber \\
    & \overset{\eqref{eqn:tilde_k}}{=} \alpha \left(\Azi - \agentexp{\Azi}\right) K(t). \label{eqn:Kt_diff}
\end{align}

\textbf{(b)}
Now, consider a minimum point $x^* \in \arg \min_{x \in \R^d} \frac{1}{m}\sum_{i = 1}^m F^i(x)$ defined by~\eqref{eqn:opt_1}. From the definition~\eqref{eqn:dist}, the local cost function $F^i(x)$ is convex for each agent $i \in \{1,\ldots,m\}$. Thus, the aggregate cost function $\sum_{i = 1}^m F^i(x)$ is also convex. Therefore, $x^* \in X^*$ if and only if~\cite{boyd2004convex}
\begin{align*}
    \nabla \sum_{i = 1}^m F^i(x^*) = 0_d,
\end{align*}
where $0_d$ denotes the $d$-dimensional zero vector. As $\nabla \sum_{i = 1}^m F^i(x) = A^T (A x - B)$, $x^*$ is the minimum point  if and only if it satisfies
\begin{align}
    A^T (A x^* - B) = 0_d. \label{eqn:new_X*}
\end{align}
From the definition of the each agent's stochastic gradient $g^{i_t}(t)$ in~\eqref{eqn:g_i} and the definition of individual cost function's gradient $g$ in~\eqref{eqn:g} we get
\begin{align}
    \agentexp{g^{\zeta_{t_t}}(t)} & = \agentexp{\Azi}x(t) - \agentexp{ \left(a^{\zeta_{t_t}}\right)^T b^{\zeta_{t_t}}} \nonumber \\
    & = \frac{1}{N}\left( \A x(t)-A^T B \right), \label{eqn:sum_grad_1}
\end{align}
where the last inequality follows from the definition of the {\em uniform} random variable $\zeta_{t_t}$.
Upon substituting from above and~\eqref{eqn:new_X*} in~\eqref{eqn:sum_grad_1} we get
\begin{align}
    \agentexp{g^{\zeta_{t_t}}(t)} & = \frac{1}{N}\A z(t). \label{eqn:sum_grad_2}
\end{align}
The R.H.S. above is the gradient of the objective cost $\frac{1}{m}\sum_{i=1}^m F^i$ evaluated at the current estimate $x(t)$ of~\eqref{eqn:opt_1}, which we denote by $\nabla F(t)$. Thus,
\begin{align}
    \nabla F(t) = \agentexp{g^{\zeta_{t_t}}(t)} & = \frac{1}{N} \A z(t). \label{eqn:exp_grad}
\end{align}
The above equation means that, the stochastic gradient $g^{\zeta_{t_t}}(t)$ at every iteration $t\geq 0$ is an unbiased estimate of the true gradient $\nabla F(t)$ given the current estimate $x(t)$.
% Recall from~\eqref{eqn:exp_grad},
% \begin{align*}
%     \nabla F(t) = \frac{1}{N} \A z(t).
% \end{align*}
Noting that $\norm{\A} = s_1$~\cite{horn2012matrix}, the above implies that
\begin{align}
    \norm{\nabla F(t)} \leq \frac{s_1}{N} \norm{z(t)}. \label{eqn:grad_norm}
\end{align}

\textbf{(c)}
Assumption~2, combined with~\eqref{eqn:exp_grad} and the definition~\eqref{eqn:variance}, implies that the conditional second moment of the stochastic gradients satisfies, for each iteration $t=0,1,\ldots$,
\begin{align}
    \agentexp{\norm{g^{\zeta_{t_t}}(t)}^2} \leq V_1 + V_G \norm{\nabla F(t)}^2, \label{eqn:sec_mom}
\end{align}
where $V_G = V_2 + 1$.
Since
\begin{align*}
    \agentexp{\norm{g^{\zeta_{t_t}}(t)}^2} = \frac{1}{N}\sum_{i=1}^N \norm{g^{{i}}(t)}^2,
\end{align*}
\eqref{eqn:sec_mom} implies that, for each $i \in \{1,\ldots,N\}$,
\begin{align}
    \norm{g^{{i}}(t)}^2 \leq V_1N + V_G N \norm{\nabla F(t)}^2. \label{eqn:sqrd_bd}
\end{align}
Thus, there exist two non-negative real scalar values $E_1 \geq \sqrt{V_1N}$ and $E_2 \geq \sqrt{V_{G}N}$ such that, for each $i \in \{1,\ldots,N\}$,
\begin{align}
    \norm{g^{i}(t)} \leq E_1 + E_2\norm{\nabla F(t)}. \label{eqn:first_mom}
\end{align}

\subsection{Notations}
\label{sub:notations}

For the positive valued parameters $\alpha$, $\delta$, and $\beta$, let
\begin{align}
    C_1 & =  \max_{i=1,\ldots,N} \norm{\Ai - \frac{1}{N}\A}, \label{eqn:c1} \\
    % C_2 & = \alpha ~ \expect{\norm{\Ai-\frac{1}{N}\A}} \norm{K_{\beta}}, \label{eqn:c2} \\
    % \rho & = \expect{\norm{I - \alpha \left(\Ai+\beta I \right)}}, \label{eqn:rho} \\
    % L & = \max_{i=1,\ldots,N} \norm{\Ai+\beta I}, \label{eqn:L} \\
    \mu & = \left(1 - \frac{2\alpha s_d}{N}(1-\alpha L)\right), \label{eqn:mu} \\
    % \sigma^2 & = \max_{j=1,\ldots,d} \E\left[\norm{\left(\Ai+\beta I \right) K_{\beta}~e_j  - e_j}^2 \right], \label{eqn:sigma} \\
    % C_3 & = \frac{\alpha N \sigma^2}{s_d \left(1-\alpha L\right)}, \label{eqn:c3} \\
    \varrho & = \norm{I - \alpha \left(\frac{1}{N}\A+\beta I \right)}, \label{eqn:varrho} \\
    C_4(t) & = (V_2 + 1)\frac{s_1^2}{N} \left(d C_3 + \norm{K_{\beta}}^2 + 2 C_2 \norm{K_{\beta}}\sum_{j=0}^t\rho^j + \norm{\widetilde{K}(0)}_F^2 \mu^{t+1} + 2 \norm{K_{\beta}} \norm{\widetilde{K}(0)}\rho^{t+1} \right), \label{eqn:c4} \\
    C_5(t) & = 2C_1 E_2 \frac{s_1}{N} \left(\norm{K_{\beta}} +  \norm{\widetilde{K}(0)} \varrho^{t} \right), \label{eqn:c5} \\
    C_6(t) & = \frac{2s_d}{s_d + N\beta} - 2 \frac{s_1}{N}\norm{\widetilde{K}(0)} \varrho^{t+1}, \label{eqn:c6} \\
    C_7(t) & = 2 C_1 E_1 \left(\norm{K_{\beta}} + \norm{\widetilde{K}(0)} \varrho^{t}\right), \label{eqn:c7} \\
    C_8(t) & = C_4(t) + 0.5, \label{eqn:c8} \\
    R_3(t) & = \delta^2 V_1 N \left(d C_3 + \norm{K_{\beta}}^2 + 2 C_2 \norm{K_{\beta}}\sum_{j=0}^t\rho^j + \norm{\widetilde{K}(0)}_F^2 \mu^{t+1} + 2 \norm{K_{\beta}} \norm{\widetilde{K}(0)}\rho^{t+1} \right), \label{eqn:r3} \\
    R_2(t) & = R_3(t) + \frac{1}{2} \alpha^2 C_7(t)^2, \label{eqn:r2} \\
    \overline{\delta}(t) & = \min \{\frac{1}{C_6(t)}, \frac{C_6(t) - \alpha C_5(t)}{C_8(t)}\}. \label{eqn:delta_bar}
\end{align}

% \subsection{Preliminary Observation}
% \label{prf:prelim}

% We have the following preliminary observations about Algorithm~\ref{algo_1} for each iteration $t\geq 0$.
% \begin{itemize}
%     \item The probability that a specific agent $l\in\{1,\ldots,m\}$ computes the stochastic gradient $g^{n(l-1)+w_l^t}(t)$, based on a particular row $w_l^t$ of its local data points $(A^l, B^l)$, is $\frac{1}{n}$.
%     \item The probability that the server picks a specific agent $k^t\in\{1,\ldots,m\}$ is $\frac{1}{m}$.
%     \item The probability that the server updates the estimate $x(t)$ based on a particular row $i_t \in \{1,\ldots,N\}$ of the collective data points $(A,B)$ is $(\frac{1}{n})(\frac{1}{m}) = \frac{1}{N}$, which is the same as selecting a data point {\em uniformly} at random from the collective data points $(A,B)$.
% \end{itemize}
% Following the above argument, the update equations for the pre-conditioner matrix in~\eqref{eqn:kcol_update_dist} and for the estimate in~\eqref{eqn:x_update_dist} are respectively equivalent to the following equations:
% \begin{align}
%     k_j(t + 1) & = k_j(t) - \alpha ~ R^{i_t}_j(t), \text{ and } \label{eqn:kcol_update} \\
%     x(t+1) & = x(t) - \delta ~ K(t+1) ~ g^{i_t}(t), \label{eqn:x_update}
% \end{align}
% where $\{i_t\}$ are drawn i.i.d. from the uniform distribution within $\{1,\ldots,N\}$, for each iteration $t=0,1,\ldots$.

\subsection{Convergence of the Pre-Conditioner Matrix}
\label{sub:lemma}

We present below a result regarding the convergence of the iterative pre-conditioner matrix in Algorithm~\ref{algo_1}.\\

To present the convergence result of the pre-conditioner matrix $K(t)$, we recall some notation from Section~\ref{sub:theorem}.
\begin{itemize}
    \item Recall that $a^i \in \R^{1\times d}$ and $b^i\in \R$ respectively denote each row $i \in \{1,\ldots,N\}$ of the collective input matrix $A$ and the collective output vector $B$.
    \item  For the positive valued
parameters $\alpha$ and $\beta$, recall from~\eqref{eqn:c2} that
\begin{align*}
    C_2 & =  \alpha ~ \frac{1}{N} \sum_{i=1}^N \norm{\Ai-\frac{1}{N}\A} \norm{K_{\beta}}.
\end{align*}
    \item For positive values of the parameters $\alpha$ and $\beta$, recall from~\eqref{eqn:rho} that
    \begin{align*}
    \rho & = \frac{1}{N} \sum_{i=1}^N \norm{I - \alpha \left(\Ai+\beta I \right)}.
\end{align*}
    \item For each iteration $ t \geq 0$, let
\begin{align}
    \widetilde{K}(t) = K(t) - K_{\beta}. \label{eqn:tilde_k}
\end{align}
    \item Recall that for each $i \in \{1,\ldots,N\}$, $\Lambda_{i}$ and $\lambda_{i}$ respectively denote the largest and the smallest eigenvalue of the positive semi-definite matrix $\Ai$. Thus, $\Lambda_i,\lambda_i > 0$.
\end{itemize}

\begin{lemma}
\label{lem:kt_norm_conv}
Consider Algorithm~\ref{algo_1} with parameter $\beta > 0$. For each iteration $t\geq 0$, if $0 < \alpha < \min_{i=1,\ldots,N} \left\{\frac{2}{\Lambda_{i} + \beta} \right\}$ then $\rho < 1$ and 
\begin{align}
    \totexp{\norm{\widetilde{K}(t+1)}} & \leq \rho^{t+1} \norm{\widetilde{K}(0)} + C_2 \sum_{j=0}^t \rho^j. \label{eqn:kt_norm_conv}
\end{align}
\end{lemma}
~\\
The proof of Lemma~\ref{lem:kt_norm_conv} is deferred to Appendix~\ref{sub:lemma}.\\

Lemma~\ref{lem:kt_norm_conv} implies that, for sufficiently small value of the parameter $\alpha$ at every iteration, the iterative pre-conditioner matrix $K(t)$ in Algorithm~\ref{algo_1} converges {\em linearly} in expectation to a neighborhood of the matrix $K_{\beta}$. Since $\rho \in (0,1)$, this neighborhood is characterized from~\eqref{eqn:kt_norm_conv} as
\begin{align*}
    \lim_{t \to \infty} \totexp{\norm{\widetilde{K}(t)}} & \leq \frac{C_2}{1-\rho}.
\end{align*}
From~\eqref{eqn:kt_norm_conv}, smaller value of the parameter $\rho$ implies faster convergence of the pre-conditioner matrix. However, the final error in $K(t)$ is large if $\rho$ is small. Thus, there is a trade-off between the rate of convergence and the final error regarding the convergence of the pre-conditioner matrix in Algorithm~\ref{algo_1}.
From~\eqref{eqn:c2}, in the deterministic case, the value of $C_2 = 0$, which means that the sequence of pre-conditioner matrices $\{K(t), ~ t=0,1,\ldots\}$ converges exactly to $K_{\beta}$ in this case.\\

\subsection{Proof of Lemma~\ref{lem:kt_norm_conv}}
\label{prf:kt_norm_conv}

Here we present the proof of Lemma~\ref{lem:kt_norm_conv}. Consider an arbitrary iteration $t\geq 0$.
\begin{itemize}
    \item Recall from~\eqref{eqn:It} that
    \begin{align*}
        I_t = \{1_t, \ldots m_t\} \cup \{\zeta_t\}. 
    \end{align*}
    \item Recall from Section~\ref{sub:assump} that for each iteration $t\geq 0$ and agent $i \in \{1,\ldots,m\}$, $\condexp{\cdot}$ denotes the conditional expectation of a function the random variable $i_t$ given the current estimate $x(t)$ and the current pre-conditioner $K(t)$. Similarly, for each iteration $t\geq 0$, $\serverexp{\cdot}$ denotes the conditional expectation of a function the random variable $\zeta_t$ given the current estimate $x(t)$ and the current pre-conditioner $K(t)$.
    Recall from~\eqref{eqn:agentexp_def} that
    \begin{align*}
        \agentexp{\cdot} = \E_{1_t, \ldots m_t,\zeta_t}(\cdot).
    \end{align*}
    \item Recall from~\eqref{eqn:def_totexp} that $\totexp{\cdot}$ denotes the total expectation of a function of the collection of the random variables $\{I_0,\ldots,I_t\}$ given the initial estimate $x(0)$ and initial pre-conditioner matrix $K(0)$. Specifically,
    \begin{align*}
        \totexp{\cdot} = \E_{I_0,\ldots,I_t}(\cdot).
    \end{align*}
    % \item For each iteration $t\geq 0$, the random variable $\zeta_{t_t}$ will be represented by the abridged notation $\zeta_t$ for the rest of this paper.
\end{itemize}
From Step 2 and Step 4 of Algorithm~\ref{algo_1}, the random variable $\zeta_{t_t}$ is {\em uniformly} distributed in the set $\{1,\ldots,N\}$ that denotes the total number of data points in $(A,B)$. Moreover, $(a^{\zeta_{t_t}},b^{\zeta_{t_t}})$ is a data point {\em uniformly} and independently drawn at random from $(A,B)$.

For each iteration $t=0,1,\ldots$, upon substituting from~\eqref{eqn:Rij} and~\eqref{eqn:R} in~\eqref{eqn:kcol_update_dist} we have, for each column index $j=1,\ldots,d$ of the matrix $K(t)$,
\begin{align}
    k_j(t + 1) & = k_j(t) - \alpha \left( \left(\Azi+ \beta I \right) k_j(t) - e_j\right). \label{eqn:kj_iter}
\end{align}
Recall the definition of $K_{\beta}$ in~\eqref{eqn:kbeta}. Let, $k_{j\beta}$ denote the $j$-th column of $K_{\beta}$. Then for each column $j=1,\ldots,d$ of $K_{\beta}$ we obtain that 
\begin{align}
    \left(\frac{1}{N}\A+\beta I\right) k_{j\beta} = e_j. \label{eqn:ej}
\end{align}
Recall the definition~\eqref{eqn:tilde_k},
\begin{align*}
    \widetilde{K}(t) = K(t) - K_{\beta}, \quad \forall t \geq 0.
\end{align*}
Let $\widetilde{k}_j(t)$ denote the $j$-th column of $\widetilde{K}(t)$. Subtracting $k_{j\beta}$ from both sides of~\eqref{eqn:kj_iter}, from the definition of $\widetilde{k}_j(t)$ we have for each $j=1,\ldots,d$,
\begin{align*}
    \widetilde{k}_j(t + 1) 
    & = \widetilde{k}_j(t) - \alpha \left( \Azi+\beta I \right) \left(\widetilde{k}_j(t) +  k_{j\beta}\right) + \alpha e_j \\
    &\overset{\eqref{eqn:ej}}{=} \widetilde{k}_j(t) - \alpha \left(\Azi+\beta I \right) \left(\widetilde{k}_j(t) +  k_{j\beta}\right)
    + \alpha \left(\frac{1}{N}\A+\beta I\right) k_{j\beta} \\
    & = \left(I - \alpha \left(\Azi+\beta I \right)\right)\widetilde{k}_j(t)
    - \alpha \left(\Azi-\frac{1}{N}\A \right)k_{j\beta}.
\end{align*}
Upon horizontally concatenating the columns $\{\widetilde{k}_j(t), ~ j=1,\ldots,d\}$, from above we get
\begin{align}
    \widetilde{K}(t+1) = \left(I - \alpha \left(\Azi+\beta I \right)\right)\widetilde{K}(t) - \alpha \left(\Azi-\frac{1}{N}\A \right)K_{\beta}. \label{eqn:Kt_iter}
\end{align}
Using triangle inequality on the R.H.S. of~\eqref{eqn:Kt_iter} we get
\begin{align*}
    \norm{\widetilde{K}(t+1)} & \leq \norm{\left(I - \alpha \left(\Azi+\beta I \right)\right)\widetilde{K}(t)} + \alpha \norm{\left(\Azi-\frac{1}{N}\A \right)K_{\beta}}.
\end{align*}
From the definition of induced 2-norm of matrix~\cite{horn2012matrix},
\begin{align*}
    \norm{\widetilde{K}(t+1)} & \leq \norm{I - \alpha \left(\Azi+\beta I \right)} \norm{\widetilde{K}(t)} + \alpha \norm{\Azi-\frac{1}{N}\A}\norm{K_{\beta}}.
\end{align*}
Upon taking conditional expectation $\agentexp{\cdot}$ (see definition~\eqref{eqn:agentexp_def}) on both sides, given the current matrix $K(t)$ and the estimate $x(t)$, we get
\begin{dmath}
    \agentexp{\norm{\widetilde{K}(t+1)}}  \leq \agentexp{\norm{I - \alpha \left(\Azi+\beta I \right)}} \norm{\widetilde{K}(t)} + \alpha \agentexp{\norm{\Azi-\frac{1}{N}\A}} \norm{K_{\beta}}. \label{eqn:kt_norm_conv_pre_1}
\end{dmath}
Since the random variable $\zeta_{t_t}$ is {\em uniformly} distributed in $\{1,\ldots,N\}$, we have
\begin{align*}
    \agentexp{\norm{I - \alpha \left(\Azi+\beta I \right)}} = \frac{1}{N} \sum_{i=1}^N \norm{I - \alpha \left(\Ai+\beta I \right)}, \\
    \agentexp{\norm{\Azi-\frac{1}{n}\A}} = \frac{1}{N} \sum_{i=1}^N \norm{\Ai-\frac{1}{N}\A}.
\end{align*}
Upon substituting from above in~\eqref{eqn:kt_norm_conv_pre_1} we obtain that
\begin{dmath}
    \agentexp{\norm{\widetilde{K}(t+1)}}  \leq \frac{1}{N} \sum_{i=1}^N \norm{I - \alpha \left(\Ai+\beta I \right)} \norm{\widetilde{K}(t)} + \alpha \frac{1}{N} \sum_{i=1}^N \norm{\Ai-\frac{1}{N}\A} \norm{K_{\beta}}. \label{eqn:kt_norm_conv_pre}
\end{dmath}
Recall from Section~\ref{sub:lemma} that $\Lambda_{i}$ and $\lambda_{i}$ respectively denote the largest and the smallest eigenvalue of each $\Ai$.
Since $\left(\Ai+\beta I \right)$ is positive definite for $\beta > 0$,
for each value of $\alpha$ satisfying $0 < \alpha < \min_{i=1,\ldots,N} \left\{\frac{2}{\Lambda_{i} + \beta} \right\}$ we have~\cite{fessler2008image}
$$\norm{I - \alpha \left(\Ai+\beta I \right)} = \max \left\{\mnorm{1-\alpha \left(\Lambda_{i} +\beta \right)}, \mnorm{1-\alpha \left(\lambda_{i} +\beta \right)}\right\} < 1, ~ i=1,\ldots,N.$$
Using the definitions of $C_2$ (see~\eqref{eqn:c2}) and $\rho$ (see~\eqref{eqn:rho}) in~\eqref{eqn:kt_norm_conv_pre} we then have $\rho < 1$ such that
\begin{align}
    \agentexp{\norm{\widetilde{K}(t+1)}} & \leq \rho \norm{\widetilde{K}(t)} + C_2. \label{eqn:kt_cond}
\end{align}
Iterating the above from $t$ to $0$, by the law of total expectation we have
\begin{align}
    \totexp{\norm{\widetilde{K}(t+1)}} & \leq \rho^{t+1} \norm{\widetilde{K}(0)} + C_2 \sum_{j=0}^t \rho^j. 
\end{align}
Hence, the proof.

%%%%%%%%%%%%%%%%%%%%%%%%%%%%%%%%%%%%%%%%%%%%%%%%%%%%%%%%%%%%%%%%%

\subsection{Proof of Theorem~\ref{thm:z_conv}}
\label{prf:z_conv}

Here, we formally prove Theorem~\ref{thm:z_conv}.

\subsubsection{Proof of Part (ii) of Theorem~\ref{thm:z_conv}}

{\bf Step I:} Upon subtracting $x^*$ from both sides of~\eqref{eqn:x_update_dist} and using the definition of $z(t)$ in~\eqref{eqn:err}, we have
\begin{align*}
    z(t+1) = z(t) - \delta ~ K(t + 1) ~ g^{\zeta_{t_t}}(t).
\end{align*}
The above implies that
\begin{align*}
    \norm{z(t+1)}^2 = & \norm{z(t)}^2 + \delta^2 \norm{K(t + 1) g^{\zeta_{t_t}}(t)}^2
    - 2\delta z(t)^T K(t + 1) g^{\zeta_{t_t}}(t).
\end{align*}
Upon taking conditional expectation $\agentexp{\cdot}$ on both sides above, given the current pre-conditioner matrix $K(t)$ and estimate $x(t)$, we have
\begin{align}
    \agentexp{\norm{z(t+1)}^2} = & \norm{z(t)}^2 + \delta^2 \agentexp{\norm{K(t + 1) g^{\zeta_{t_t}}(t)}^2} - 2\delta z(t)^T \agentexp{K(t + 1) g^{\zeta_{t_t}}(t)}. \label{eqn:zt_total_1}
\end{align}
Consider the expression $z(t)^T \agentexp{K(t + 1) g^{\zeta_{t_t}}(t)}$ above. Upon substituting from~\eqref{eqn:tilde_k} we have
\begin{align}
    z(t)^T \agentexp{K(t + 1) g^{\zeta_{t_t}}(t)}
    & = z(t)^T \agentexp{\widetilde{K}(t + 1) g^{\zeta_{t_t}}(t)} + z(t)^T K_{\beta}\agentexp{g^{\zeta_{t_t}}(t)} \nonumber\\
    & \overset{\eqref{eqn:exp_grad}}{=} z(t)^T \agentexp{\widetilde{K}(t + 1) g^{\zeta_{t_t}}(t)} + z(t)^T K_{\beta}\nabla F(t). \label{eqn:zt_1}
\end{align}
Consider the following equation
\begin{align*}
    z(t)^T \agentexp{\agentexp{\widetilde{K}(t + 1)} \left(\nabla F(t) - g^{\zeta_{t_t}}(t)\right)}
    & = z(t)^T \agentexp{\widetilde{K}(t + 1)} \agentexp{\left(\nabla F(t) - g^{\zeta_{t_t}}(t)\right)} \overset{\eqref{eqn:exp_grad}}{=} 0.
\end{align*}
Owing to the above, we have the first expression in the R.H.S. of~\eqref{eqn:zt_1} as
\begin{align}
    & z(t)^T \agentexp{\widetilde{K}(t + 1) g^{\zeta_{t_t}}(t)} \nonumber \\
    = & z(t)^T \agentexp{\widetilde{K}(t + 1)} \nabla F(t) - z(t)^T \agentexp{\widetilde{K}(t + 1) \left(\nabla F(t) - g^{\zeta_{t_t}}(t)\right)} \nonumber \\
    & + z(t)^T \agentexp{\agentexp{\widetilde{K}(t + 1)} \left(\nabla F(t) - g^{\zeta_{t_t}}(t)\right)} \nonumber\\
    = & z(t)^T \agentexp{\widetilde{K}(t + 1)} \nabla F(t) - z(t)^T \agentexp{\widetilde{K}(t + 1)-\agentexp{\widetilde{K}(t + 1)}}\nabla F(t) \nonumber\\
    & - z(t)^T \agentexp{\left(\agentexp{\widetilde{K}(t + 1)} - \widetilde{K}(t + 1)\right) g^{\zeta_{t_t}}(t)} \nonumber \\
    = & z(t)^T \agentexp{\widetilde{K}(t + 1)} \nabla F(t) - z(t)^T \agentexp{\left(\agentexp{\widetilde{K}(t + 1)} - \widetilde{K}(t + 1)\right) g^{\zeta_{t_t}}(t)}. \label{eqn:zt_2}
\end{align}
From~\eqref{eqn:Kt_diff} we have the second expression in the R.H.S. above as
\begin{dmath*}
    z(t)^T \agentexp{\left(\agentexp{\widetilde{K}(t + 1)} - \widetilde{K}(t + 1)\right) g^{\zeta_{t_t}}(t)} = \alpha z(t)^T \agentexp{\left(\Azi - \agentexp{\Azi}\right)K(t) g^{\zeta_{t_t}}(t)}.
\end{dmath*}
Applying Cauchy-Schwartz inequality~\cite{roman2005advanced} above we have
\begin{dmath*}
    z(t)^T \agentexp{\left(\agentexp{\widetilde{K}(t + 1)} - \widetilde{K}(t + 1)\right) g^{\zeta_{t_t}}(t)}  \leq \alpha \norm{z(t)} \norm{ \agentexp{\left(\Azi - \agentexp{\Azi}\right)K(t)g^{\zeta_{t_t}}(t)}}.
\end{dmath*}
Using Jensen's inequality~\cite{rudin2006real} on the convex function $\norm{\cdot}$, from above we have
\begin{dmath}
    z(t)^T \agentexp{\left(\agentexp{\widetilde{K}(t + 1)} - \widetilde{K}(t + 1)\right) g^{\zeta_{t_t}}(t)} \leq \alpha \norm{z(t)} \agentexp{\norm{\left(\Azi - \agentexp{\Azi}\right)K(t) g^{\zeta_{t_t}}(t)}}. \label{eqn:jensen}
\end{dmath}
From the definition of induced $2$-norm of matrix we have
\begin{dmath*}
    \norm{\left(\Azi - \agentexp{\Azi}\right)K(t) g^{\zeta_{t_t}}(t)} \leq \norm{\Azi - \agentexp{\Azi}} \norm{K(t)} \norm{g^{\zeta_{t_t}}(t)} \overset{\eqref{eqn:c1}}{\leq} C_1 \norm{K(t)} \norm{g^{\zeta_{t_t}}(t)}.
\end{dmath*}
From the definition of expectation and substituting from above in~\eqref{eqn:jensen} we get
\begin{align}
    z(t)^T \agentexp{\left(\agentexp{\widetilde{K}(t + 1)} - \widetilde{K}(t + 1)\right) g^{\zeta_{t_t}}(t)}
    & \leq \alpha C_1 \norm{z(t)} \norm{K(t)} \agentexp{\norm{g^{\zeta_{t_t}}(t)}} \nonumber\\
    & \overset{\eqref{eqn:first_mom}}{\leq} \alpha C_1 \norm{z(t)} \norm{K(t)} \left( E_1 + E_2\norm{\nabla F(t)}\right).  \label{eqn:zt_3}
\end{align}
Upon substituting from~\eqref{eqn:zt_1},~\eqref{eqn:zt_2} and~\eqref{eqn:zt_3} in~\eqref{eqn:zt_total_1} we obtain that
\begin{dmath}
    \agentexp{\norm{z(t+1)}^2}
    \leq \norm{z(t)}^2 + \delta^2 \agentexp{\norm{K(t + 1) g^{\zeta_{t_t}}(t)}^2} 
    - 2\delta z(t)^T K_{\beta}\nabla F(t) - 2\delta z(t)^T \agentexp{\widetilde{K}(t + 1)} \nabla F(t) 
    + 2\delta \alpha C_1 \norm{z(t)} \norm{K(t)} \left( E_1 + E_2\norm{\nabla F(t)}\right). \label{eqn:zt_total_2}
\end{dmath}
In the following Steps II-V, we bound the expressions in the R.H.S. above.\\

{\bf Step II:} In this step, we bound the second expression $\delta^2 \agentexp{\norm{K(t + 1) g^{\zeta_{t_t}}(t)}^2}$ in the R.H.S. of~\eqref{eqn:zt_total_2}.\\

Consider the expression $\agentexp{\norm{K(t + 1) g^{\zeta_{t_t}}(t)}^2}$. Using the definition of induced $2$-norm and expectation we have
\begin{align}
    \agentexp{\norm{K(t + 1) g^{\zeta_{t_t}}(t)}^2} & \leq \agentexp{\norm{K(t+1)}^2 \norm{g^{\zeta_{t_t}}(t)}^2} \nonumber \\
    & \overset{\eqref{eqn:sqrd_bd}}{\leq} \left( V_1 N + V_G N \norm{\nabla F(t)}^2 \right) \agentexp{\norm{K(t+1)}^2}. \label{eqn:kt_gi_total}
\end{align}
Using triangle inequality on induced $2$-norm in~\eqref{eqn:tilde_k} we get
\begin{align}
    \norm{K(t)} \leq \norm{\widetilde{K}(t)} + \norm{K_{\beta}}, \, \forall t\geq 0. \label{eqn:norm_ktilde}
\end{align}
From~\eqref{eqn:norm_ktilde} and the definition of expectation,
\begin{align}
    \agentexp{\norm{K(t+1)}^2} & \leq  \agentexp{\norm{\widetilde{K}(t+1)}^2} + \norm{K_{\beta}}^2 + 2 \norm{K_{\beta}} \agentexp{\norm{\widetilde{K}(t+1)}}. \label{eqn:kt_gi_1}
\end{align}
In the rest of this step, we bound the first expression above, which in turn bounds the R.H.S. of~\eqref{eqn:kt_gi_total}. Note that, the third expression above has already been bounded in~\eqref{eqn:kt_cond}.\\

% Using triangle inequality on the R.H.S. of~\eqref{eqn:Kt_iter} we get
% \begin{align*}
%     \norm{\widetilde{K}(t+1)} & \leq \norm{I - \alpha \left(\Ai+\beta I \right)} \norm{\widetilde{K}(t)} + \alpha \norm{\Ai-\frac{1}{N}\A}\norm{K_{\beta}}.
% \end{align*}
% Then we have
% \begin{align}
%     \condexp{\norm{\widetilde{K}(t+1)}} & \leq \expect{\norm{I - \alpha \left(\Ai+\beta I \right)}} \norm{\widetilde{K}(t)} + \alpha \expect{\norm{\Ai-\frac{1}{N}\A}} \norm{K_{\beta}}. \label{eqn:kt_norm_conv_pre}
% \end{align}
% Since $\left(\Ai+\beta I \right)$ is positive definite for $\beta > 0$,
% for each value of $\alpha$ satisfying $0 < \alpha < \max_{i = 1,\ldots,N} \left\{\frac{2}{\Lambda_i + \beta} \right\}$ we have
% $\norm{I - \alpha \left(\Ai+\beta I \right)} = \max \left\{\mnorm{1-\alpha \left(\Lambda_i +\beta \right)}, \mnorm{1-\alpha \left(\lambda_i +\beta \right)}\right\} < 1$ (Fessler, 2020). Using~\eqref{eqn:c2} and the definition of $\rho$ (see~\eqref{eqn:rho}) in~\eqref{eqn:kt_norm_conv_pre} we then have $\rho < 1$ such that
% \begin{align*}
%     \condexp{\norm{\widetilde{K}(t+1)}} & \leq \rho \norm{\widetilde{K}(t)} + C_2. 
% \end{align*}
% Iterating the above from $t$ to $0$,
% \begin{align}
%     \condexp{\norm{\widetilde{K}(t+1)}} & \leq \rho^{t+1} \norm{\widetilde{K}(0)} + C_2 \sum_{j=0}^t \rho^j. \label{eqn:kt_norm_conv}
% \end{align}

For each $i \in \{1,\ldots,N\}$ and each $j \in \{1,\ldots,d\}$,
define a function $h^i_j: \R^d \to \R $ such that
\begin{align}
    h^i_j(x) = \frac{1}{2} x^T \left(\Ai+\beta I \right) x - x^Te_j. \label{eqn:hij}
\end{align}
The gradient of $h^i_j$ is given by
\begin{align}
    \nabla h^i_j(x) = \left(\Ai+\beta I \right) x - e_j, \, \forall x\in \R^d. \label{eqn:hij_grad}
\end{align}
Note that,
\begin{align}
    \agentexp{\nabla h^{\zeta_{t_t}}_j(x)} & = \left(\agentexp{\Azi}+ \beta I \right) x - e_j  = \left(\frac{1}{N}\A + \beta I \right) x - e_j, \, \forall x\in \R^d. \label{eqn:hj_grad}
\end{align}
Hence,~\eqref{eqn:kj_iter} is a stochastic gradient descent update with stepsize $\alpha$ on the objective cost function $H_j: \R^d \to \R $ defined by $H_j(x) = \frac{1}{2}x^T \left(\frac{1}{N}\A + \beta I \right) x - x^T e_j$. Under Assumption 1, the matrix $\A$ is symmetric positive definite. Equivalently, $H_j$ is $\frac{s_d}{N}$-strongly convex. Also note that, each $\nabla h^i_j$ is Lipschitz with a Lipschitz constant $L_i = \norm{\Ai+\beta I} = \Lambda_i + \beta$.\\

Upon substituting from~\eqref{eqn:ej} into~\eqref{eqn:hj_grad} we have $\agentexp{\nabla h^{\zeta_{t_t}}_j(k_{j\beta})} = 0_d$, where $0_d$ denotes the origin of $\R^d$. For each $j=1,\ldots,d$, we define a quantity
\begin{align}
    \sigma_j^2 & = \max_{t\geq 0} \agentexp{\norm{\nabla h^{\zeta_{t_t}}_j(k_{j\beta})}^2} \label{eqn:sigma_j}
\end{align}
Upon substituting above from~\eqref{eqn:hij_grad} and~\eqref{eqn:ej} we get
\begin{align*}
    \sigma_j^2 & = \max_{t\geq 0} \agentexp{\norm{\left(\Azi+\beta I \right) K_{\beta}~ e_j  - e_j}^2} \\
    & = \frac{1}{N} \sum_{i=1}^N \norm{\left(\Ai+\beta I \right) K_{\beta}~ e_j  - e_j}^2.
\end{align*}
For each $j=1,\ldots,d$, let $\widetilde{k}_j(t) = k_j(t) - k_{j\beta}$. 
Then, if $\alpha < \min \left\{\frac{N}{s_d},\frac{1}{L}\right\}$, we have for each $j=1,\ldots,d$,~\cite{needell2014stochastic}
\begin{align}
    \E_t \left[\norm{\widetilde{k}_j(t+1)}^2 \right] \leq \mu^{t+1} \norm{\widetilde{k}_j(0)}^2 + C_3. \label{eqn:k_conv_pre}
\end{align}
where $\mu \in (0,1)$ (see~\eqref{eqn:mu}).
Let $\norm{\cdot}_F$ denote the Frobenius norm of a matrix~\cite{horn2012matrix}. Then,
\begin{align*}
    \norm{\widetilde{K}(t)}^2 \leq \norm{\widetilde{K}(t)}^2_F = \sum_{j=1}^d \norm{\widetilde{k}_j(t)}^2,
\end{align*}
which implies that
\begin{align*}
    \agentexp{\norm{\widetilde{K}(t)}^2} \leq \sum_{j=1}^d \agentexp{\norm{\widetilde{k}_j(t)}^2}.
\end{align*}
Taking the total expectation above and substituting from~\eqref{eqn:k_conv_pre} we obtain that
\begin{align}
    \totexp{\norm{\widetilde{K}(t+1)}^2} \leq  \mu^{t+1}\norm{\widetilde{K}(0)}_F^2 + d C_3. \label{eqn:k_conv}
\end{align}
Upon substituting from~\eqref{eqn:kt_norm_conv} and~\eqref{eqn:k_conv} in~\eqref{eqn:kt_gi_1} we have
\begin{align*}
    & \totexp{\norm{K(t+1)}^2} \leq  \mu^{t+1}\norm{\widetilde{K}(0)}_F^2 + d C_3 + \norm{K_{\beta}}^2 + 2 \norm{K_{\beta}} \left( \rho^{t+1} \norm{\widetilde{K}(0)} + C_2 \sum_{j=0}^t \rho^j \right).
\end{align*}
Upon substituting from above in~\eqref{eqn:kt_gi_total} we obtain that
\begin{dmath}
    \totexp{\norm{K(t + 1) g^{i_t}(t)}^2} 
    \leq \left( V_1N + V_G N \norm{\nabla F(t)}^2 \right) \left(\mu^{t+1}\norm{\widetilde{K}(0)}_F^2 + d C_3
    + \norm{K_{\beta}}^2 + 2 \norm{K_{\beta}} \left( \rho^{t+1} \norm{\widetilde{K}(0)} + C_2 \sum_{j=0}^t \rho^j \right)\right). \label{eqn:kt_gi}
\end{dmath}

{\bf Step III:} In this step, we bound the the third expression $-z(t)^T K_{\beta}\nabla F(t)$ in the R.H.S. of~\eqref{eqn:zt_total_2}.\\

From eigen value decomposition~\cite{horn2012matrix}, 
$\A = V S V^T$
where 
$S = Diag\left(s_1, \ldots,  s_d \right)$,
and the matrix $V$ constitutes of orthonormal eigen vectors of $\A$. From above, $$\left(\frac{1}{N}\A+\beta I\right)^{-1} \A  = V Diag\left( \frac{s_1}{(s_1/N) + \beta}, \ldots, \frac{s_d}{(s_d/N) + \beta}\right) V^T.$$ 
Since $K_{\beta} = \left(\frac{1}{N}\A+\beta I\right)^{-1}$ (see~\eqref{eqn:kbeta}), the eigenvalues of $K_{\beta}\A$ are given by $\left \{ \frac{s_i}{(s_i/N) + \beta} | \, i=1,\ldots,d \right\}$. Thus, from the bound on Rayleigh quotient~\cite{horn2012matrix}
\begin{align}
    - z(t)^T K_{\beta}\A z(t) \leq -\lambda_{min}(K_{\beta}\A) \norm{z(t)}^2, \label{eqn:neg_term_1}
\end{align}
where $\lambda_{min}(K_{\beta}\A) = \frac{s_d}{(s_d/N) + \beta}$ is the minimum eigen value of $K_{\beta}\A$.
Upon substituting from~\eqref{eqn:exp_grad} in the expression $-z(t)^T K_{\beta}\nabla F(t)$, we obtain that
\begin{align}
    -z(t)^T K_{\beta}\nabla F(t) = - \frac{1}{N} z(t)^T K_{\beta}\A z(t) \overset{\eqref{eqn:neg_term_1}}{\leq} -\frac{1}{N}\lambda_{min}(K_{\beta}\A) \norm{z(t)}^2. \label{eqn:neg_term}
\end{align}

{\bf Step IV:} In this step, we bound $\norm{\widetilde{K}(t)}$ which appears in the fifth expression in the R.H.S. of~\eqref{eqn:zt_total_2}.\\

Recall the definition of the {\em uniform} random variable $\zeta_{t_t}$. Then,
\begin{align*}
    \agentexp{\Azi} = \frac{1}{N}\A.
\end{align*}
Upon substituting from above in~\eqref{eqn:Kt_iter_exp} we have
\begin{align*}
    \agentexp{\widetilde{K}(t+1)}
    & = \left(I - \alpha \left(\agentexp{\Azi}+\beta I \right)\right)\widetilde{K}(t) = \left(I - \alpha \left(\frac{1}{N}\A +\beta I \right)\right)\widetilde{K}(t).
\end{align*}
Again using the definition of induced matrix $2$-norm above, we obtain that
\begin{align}
    \norm{\agentexp{\widetilde{K}(t+1)}} & \leq \norm{I - \alpha \left(\frac{1}{N}\A+\beta I \right)} \norm{\widetilde{K}(t)}. \label{eqn:kt_iter_2}
\end{align}
Since $\left(\frac{1}{N}\A+\beta I \right)$ is positive definite for $\beta > 0$, for each value of $\alpha$ satisfying $0 < \alpha < \frac{2}{s_1/N + \beta}$  we have
$\norm{I - \alpha \left(\frac{1}{N}\A+\beta I \right)} = \max \left\{\mnorm{1-\alpha \left(\frac{s_1}{N}+\beta \right)}, \mnorm{1-\alpha \left(\frac{s_d}{N}+\beta \right)}\right\} < 1$ (Fessler, 2020). Using the definition of $\varrho$ (see~\eqref{eqn:varrho}) in~\eqref{eqn:kt_iter_2} we then have $\varrho < 1$ such that
\begin{align*}
    \norm{\agentexp{\widetilde{K}(t+1)}} & \leq \varrho \norm{\widetilde{K}(t)}.
\end{align*}
Iterating the above from $t$ to $0$, the total expectation is given by
\begin{align}
    \norm{\E_t \left[\widetilde{K}(t+1) \right]} & \leq \varrho^{t+1} \norm{\widetilde{K}(0)}. \label{eqn:kt_iter_3}
\end{align}
From the definition of conditional expectation $\Et[\cdot]$, we have $\Et \left[\widetilde{K}(t) \right] = \widetilde{K}(t)$. Thus,~\eqref{eqn:kt_iter_3} implies that
\begin{align}
    \norm{\widetilde{K}(t)} & \leq \varrho^t \norm{\widetilde{K}(0)}. \label{eqn:expected_kt}
\end{align}

{\bf Step V:} In this step, we bound the the fourth expression $-z(t)^T \agentexp{\widetilde{K}(t + 1)} \nabla F(t)$ in the R.H.S. of~\eqref{eqn:zt_total_2}.\\

Using Cauchy-Schwartz inequality~\cite{roman2005advanced},
\begin{align}
    -z(t)^T \agentexp{\widetilde{K}(t + 1)} \nabla F(t) \leq \norm{z(t)}\norm{\nabla F(t)} \norm{\agentexp{\widetilde{K}(t + 1)}} \overset{\eqref{eqn:kt_iter_3}}{\leq} \varrho^{t+1} \norm{\widetilde{K}(0)}\norm{z(t)}\norm{\nabla F(t)}. \label{eqn:cs}
\end{align}

{\bf Step VI:} In this step, we combine the upper bounds obtained in Step-II to Step-V above to get a bound on the R.H.S. of~\eqref{eqn:zt_total_2} in Step-I.\\

Upon substituting from~\eqref{eqn:exp_grad},~\eqref{eqn:norm_ktilde},~\eqref{eqn:kt_gi},~\eqref{eqn:neg_term},~\eqref{eqn:cs},~\eqref{eqn:expected_kt} in~\eqref{eqn:zt_total_2} we obtain that
\begin{align}
    \totexp{\norm{z(t+1)}^2} \leq \left( 1 + \delta^2 C_4(t) + \alpha \delta C_5(t) - \delta C_6(t) \right) \norm{z(t)}^2 + \alpha \delta C_7(t) \norm{z(t)} + R_3(t), \label{eqn:zt_4}
\end{align}
where $C_{4-7}(t)$ and $R_3(t)$ have been defined in~\eqref{eqn:c4}-\eqref{eqn:r3}. 
We apply the AM-GM inequality~\cite{rudin2006real} $2ab \leq (a^2 + b^2)$ with $a = \delta \norm{z(t)}$ and $b = \alpha C_7(t)$, which results in
\begin{align*}
    \totexp{\norm{z(t+1)}^2} \leq \left( 1 + \delta^2 (C_4(t) + 0.5) + \alpha \delta C_5(t) - \delta C_6(t) \right) \norm{z(t)}^2 + R_3(t) + \frac{1}{2} \alpha^2 C_7(t)^2.
\end{align*}
From the definitions of $C_8(t)$ in~\eqref{eqn:c8} and $R_2(t)$ in~\eqref{eqn:r2}, the above can be rewritten as
\begin{align}
    \totexp{\norm{z(t+1)}^2} & \leq \left( 1 + \delta^2 C_8(t) + \alpha \delta C_5(t) - \delta C_6(t) \right) \norm{z(t)}^2 + R_2(t)\nonumber\\
    & \overset{\eqref{eqn:r1}}{=} R_1(t) \norm{z(t)}^2 + R_2(t). \label{eqn:zt_5}
\end{align}
From~\eqref{eqn:kt_cond},~\eqref{eqn:k_conv_pre}, and~\eqref{eqn:kt_iter_3}, we have shown that~\eqref{eqn:zt_4}, and hence~\eqref{eqn:zt_5}, hold for any $\delta > 0$ and $\alpha$ satisfying
\begin{align*}
    0 & < \alpha < \min \left\{\frac{N}{s_d}, \frac{1}{L}, \frac{2}{(s_1/N)+\beta}, \frac{2}{\Lambda_i + \beta} | i = 1,\ldots,N  \right\}.
\end{align*}
By definition, $L = \max_{i=1,\ldots,N} \{\Lambda_i+\beta\} $ (see~\eqref{eqn:L}). Thus, the above implies that~\eqref{eqn:zt_5} holds for
\begin{align}
    0 & < \alpha < \min \left\{\frac{N}{s_d}, \frac{1}{L}, \frac{2}{(s_1/N)+\beta}\right\}, ~ \delta > 0. \label{eqn:alpha}
\end{align}
From~\eqref{eqn:alpha_bar}, the upper bound of $\alpha$ above is $\overline{\alpha}$. Then,~\eqref{eqn:zt_5} completes the first part of the proof.

\subsubsection{Proof of Part (i) of Theorem~\ref{thm:z_conv}}

Next, we show that the value of $R_1(t)$ is between $0$ and $1$ after some finite iteration if the parameter $\delta$ satisfies an additional criterion.

Note that, the values of $C_5(t)$, $C_6(t)$ and $C_8(t)$ implicitly depend on $\alpha$, but independent of $\delta$. For any $\alpha > 0$ consider the function $r_{\alpha}: \mathcal{Z}^+ \to \R$ defined as
\begin{align*}
    r_{\alpha}(t) = \frac{C_6(t)}{\alpha C_5(t)} & = \frac{\frac{s_d}{s_d + N\beta} -  \frac{s_1}{N}\norm{\widetilde{K}(0)} \varrho^{t+1}}{\alpha C_1 E_2 \frac{s_1}{N} \left(\norm{K_{\beta}} +  \norm{\widetilde{K}(0)} \varrho^{t} \right)}.
\end{align*}
From the definition, $r_{\alpha}$ is a strictly increasing function.
As $\rho \in (0,1)$, we have the limit of $r_{\alpha}$ as
\begin{align*}
    \lim_{t \to \infty} r_{\alpha}(t) = \frac{C_6(t)}{\alpha C_5(t)} & = \frac{\frac{s_d}{s_d + N\beta}}{\alpha C_1 E_2 \frac{s_1}{N} \norm{K_{\beta}}} > 0,
\end{align*}
where the last inequality holds under Assumption~1.
Since $r_{\alpha}(t)$ is strictly increasing with $t$ with a positive limit as $t \to \infty$, there exists $0 \leq T < \infty$ such that $r_{\alpha}(t)>0$ for all $t \geq T$. Since the above holds for any $\alpha > 0$, it also holds for the values of $\alpha$ satisfying~\eqref{eqn:alpha}. Hence, for any value of the parameter $\alpha$ in~\eqref{eqn:alpha} there exists $T < \infty$ such that
\begin{align*}
    C_6(t) - \alpha C_5(t) > 0, ~ \forall t \geq T.
\end{align*}

Now consider any iteration $t \geq T$.
If $0 < \delta < \frac{C_6(t) - \alpha C_5(t)}{C_8(t)}$, then we have
\begin{align*}
    & \delta C_8(t) + \alpha C_5(t) < C_6(t) \\
    \implies & 1 + \delta^2 C_8(t) +  \delta \alpha C_5(t) - \delta C_6(t) < 1.
\end{align*}
Additionally, if $0 < \delta < \frac{1}{C_6(t)}$, then
\begin{align*}
    & \delta C_6(t) < 1 \\
    \implies & 1 - \delta C_6(t) > 0 \\
    \implies & 1 + \delta^2 C_8(t) +  \delta \alpha C_5(t) - \delta C_6(t) > 0.
\end{align*}
We combine the above two results to conclude the proof. If $\alpha$ satisfies~\eqref{eqn:alpha} then there exists a non-negative integer $T < \infty$ such that for any iteration $t \geq T$, if $\delta < \min \{\frac{1}{C_6(t)}, \frac{C_6(t) - \alpha C_5(t)}{C_8(t)}\}$ then 
$$\left(1 + \delta^2 C_8(t) +  \delta \alpha C_5(t) - \delta C_6(t)\right) \in (0,1).$$
From~\eqref{eqn:delta_bar}, note that the bound on $\delta$ above is $\overline{\delta}$. Thus, the proof of the second part of the theorem is complete.

\subsubsection{Proof of Part (iii) of Theorem~\ref{thm:z_conv}}

From Part (ii) of Theorem~\ref{thm:z_conv}, we have $0 < R_1(t) < 1$ for all $ t\geq T$. Thus, upon retracing~\eqref{eqn:iter_err} from $t$ to $0$, we obtain the limiting value of the expected estimation error when $t \to \infty$ as
\begin{align*}
    \lim_{t \to \infty} \totexp{\norm{z(t+1}^2} & \leq R_2(\infty).
\end{align*}
Upon substituting from~\eqref{eqn:c7} in the definition of $R_2(t)$ in~\eqref{eqn:r2}, we get
\begin{dmath*}
    R_2(\infty) = \lim_{t \to \infty} \delta^2 V_1 N \left(d C_3 + \norm{K_{\beta}}^2 + 2 C_2 \norm{K_{\beta}}\sum_{j=0}^t\rho^j + \norm{\widetilde{K}(0)}_F^2 \mu^{t+1} + 2 \norm{K_{\beta}} \norm{\widetilde{K}(0)}\rho^{t+1} \right) + \frac{1}{2} \alpha^2 \left(2 C_1 E_1 \left(\norm{K_{\beta}} + \norm{\widetilde{K}(0)} \varrho^{t}\right) \right)^2.
\end{dmath*}
From the proof of Part (i) of Theorem~\ref{thm:z_conv}, the values of $\rho, \mu, \varrho$ are within the range $(0,1)$.
Then, we have
\begin{dmath*}
    R_2(\infty) = \delta^2 V_1 N \left(d C_3 + \norm{K_{\beta}}^2 + \frac{2 C_2 \norm{K_{\beta}}}{1-\rho} +  \right) + 2 \alpha^2 \left(C_1 E_1 \norm{K_{\beta}} \right)^2.
\end{dmath*}
Hence, the proof.

\newpage
\section{APPLICATION: DISTRIBUTED STATE ESTIMATION}
\label{sec:estimation}

In this section, we formulate distributed state estimation problem in linear dynamical systems~\cite{mitra2018distributed,wang2017distributed,park2016design} as a distributed linear regression problem.\\

First, we revisit the distributed linear regression below.
Consider a {\em server-based} system architecture that comprises $m$ agents and one server. The agents can only interact with the server, and the overall system is assumed synchronous. Each agent $i$ has $n_i$ {\em local} data points, represented by an {\em input matrix} $A^i$ and an {\em output vector} $b^i$ of dimensions $n_i \times d$ and $n_i \times 1$, respectively. Thus, for all $i \in \{1, \ldots, \, m\}$, $A^i \in \R^{n_i \times d}$ and $b^i\in \R^{n_i}$. 
Note that the agents \underline{do not} share their individual local data points with the server. For each agent $i$, we define a {\em local cost function} $F^i: \R^d \to \R$ such that for a given {\em regression parameter} $x \in \R^d$,
\begin{align}
     F^i(x) = \frac{1}{2} \norm{A^i x - b^i}^2, \label{eqn:loc_cost}
\end{align}
where $\norm{\cdot}$ denotes the Euclidean norm.
A distributed regression algorithm prescribes instructions for coordination of the server and the agents to compute an optimal regression parameter $x^* \in \R^d$ such that
\begin{align}
  x^* \in X^* = \arg \min_{x \in \R^d} \, \sum_{i = 1}^m F^i(x). \label{eqn:opt}
\end{align}

Next, we give a brief description of the distributed state estimation problem in linear dynamical systems. Consider a $d$-dimensional linear time-invariant (LTI) system that is governed by the following dynamics in discrete-time:
\begin{align}
    z(t+1) = \Ax z(t), ~ t \in 0,1,\ldots, \label{eqn:dyn}
\end{align}
where $z(t) \in \mathbb{R}^d$. The {\em state matrix} of the LTI system~\eqref{eqn:dyn}, denoted by $\Ax$, is a $(d \times d)$-dimensional real square matrix. Thus, $\Ax \in \R^{d \times d}$.
Each agent $i \in \{1,\ldots,m\}$ in the system has $d$ scalar {\em local measurements} of the above system, denoted by $\{\Bar{y}^i(t) \in \R, ~ t=0,\ldots,d-1\}$,
\begin{align}
    \Bar{y}^i(t) = C^i z(t), ~ t = 0,\ldots,d-1. \label{eqn:loc_out}
\end{align}
The {\em local measurement vector} of each agent $i$, denoted by $C^i$, is a $d$-dimensional real row vector. Thus, $C^i \in \R^{1\times d}$. So, each agent $i$ has a pair of information $(\Ax, C^i)$. The task of the agents is to estimate the state of the system $z(t)$, without sharing their local measurement vectors. Since an agent do not have access to the {\em collective measurements} $\{\Bar{y}^1(t),\ldots,\Bar{y}^m(t), ~ t=0,\ldots,d-1\}$, they collaborate with a server and estimate $z(t)$ in a distributed manner. The aforementioned problem is referred as {\em distributed state estimation}~\cite{mitra2018distributed,wang2017distributed,park2016design}.\\

We define the following notations.
\begin{itemize}
    \item The {\em local observability matrix} of each agent $i \in \{1,\ldots,m\}$ is defined as
    \begin{align}
        O^i = \begin{bmatrix} (C^{i})^T & (C^i \Ax)^T & \ldots & (C^i \Ax^{d-1})^T \end{bmatrix}^T. \label{eqn:loc_obsv}
    \end{align}
    Note that, $O^i \in \mathbb{R}^{d \times d}$.
    \item The {\em global measurement matrix} is defined as
    \begin{align}
        C = \begin{bmatrix} (C^{1})^T  & \ldots & (C^{m})^T \end{bmatrix}^T. \label{eqn:gl_out}
    \end{align}
    Note that, $C \in \R^{m \times d}$.
    \item The {\em global observability matrix} is defined as
    \begin{align}
        \Bar{O} = \begin{bmatrix} C^{T} & (CA)^T & \ldots (CA^{d-1})^T \end{bmatrix}^T. \label{eqn:gl_obsv}
    \end{align}
    Note that $\Bar{O} \in \R^{m d \times d}$.
    \item For each agent $i \in \{1,\ldots,m\}$ define the following column vector $y^i$ upon stacking the local measurements $\{\Bar{y}^i(t) \in \R, ~ t=0,\ldots,d-1\}$ as
    \begin{align}
        y^i = \begin{bmatrix} \Bar{y}^i(0) & \ldots & \Bar{y}^i(d-1) \end{bmatrix}^T. \label{eqn:yi_def}
    \end{align}
    Note that, $y^i \in \R^d$.
    \item Let $\Bar{Y}$ denote the column vector obtained upon stacking the local measurements of each agent as follows:
    \begin{align}
        \Bar{Y} = \begin{bmatrix} \Bar{y}^1(0)\ldots \Bar{y}^m(0)\ldots \Bar{y}^1(d-1)\ldots \Bar{y}^m(d-1) \end{bmatrix}^T. \label{eqn:bar_y_def}
    \end{align}
    Note that $\Bar{Y} \in \R^{m d}$.
\end{itemize}

\begin{assumption} \label{assump:obsv}
We assume that, the {\em local system} $(\Ax,C^i)$ for each agent $i\in \{1,\ldots,m\}$ is unobservable, and the {\em global system} $(\Ax,C)$ is jointly observable. In other words, the local observability matrix $O^i$ of each agent $i$ is rank-deficient and the global observability matrix $\Bar{O}$ is full rank.
\end{assumption}

Assumption~\ref{assump:obsv} is standard for the distributed state estimation problem of LTI systems~~\cite{mitra2018distributed,wang2017distributed,park2016design}. \\

From~\eqref{eqn:dyn}, note that the system state at sampling time $t \geq 0$ is given by
\begin{align}
    z(t) = \Ax^t z(0), \label{eqn:dyn_iter}
\end{align}
with the understanding that $\Ax^0 = I$.
Thus, an agent can estimate $z(t)$ for all sampling time $t\geq 0$ using~\eqref{eqn:dyn_iter} if the estimate of the initial state $z(0)$ is available. Based on this fact, we frame the aforementioned distributed state estimation problem as a distributed linear regression problem~\eqref{eqn:opt} as follows.\\

Upon substituting from~\eqref{eqn:gl_out} in~\eqref{eqn:gl_obsv} we have
\begin{align}
    \Bar{O} & = \begin{bmatrix}
    (C^{1})^T\ldots (C^{m})^T\ldots (C^1\Ax^{d-1})^T\ldots (C^m\Ax^{d-1})^T \end{bmatrix}^T. \label{eqn:obar_expand}
\end{align}
Upon substituting from~\eqref{eqn:loc_out} and~\eqref{eqn:dyn_iter} in~\eqref{eqn:bar_y_def} we obtain that
\begin{align}
  \Bar{Y} & = \begin{bmatrix}
    (C^{1})^T\ldots (C^{m})^T\ldots (C^1\Ax^{d-1})^T\ldots (C^m\Ax^{d-1})^T \end{bmatrix}^T z(0) \nonumber \\
    & \overset{\eqref{eqn:obar_expand}}{=} \Bar{O}z(0). \label{eqn:bar_y}
\end{align}
We define the following vector $Y$, which is a rearrangement of the rows in $\Bar{Y}$ (ref.~\eqref{eqn:bar_y_def}) in so that the measurements of each agent $i \in \{1,\ldots,m\}$ are stacked consecutively as follows:
\begin{align}
   Y & = \begin{bmatrix} \Bar{y}^1(0)\ldots \Bar{y}^1(d-1)\ldots \Bar{y}^m(0)\ldots \Bar{y}^m(d-1) \end{bmatrix}^T. \label{eqn:y_def}
\end{align}
Similarly, we define the following matrix $O$, which is a rearrangement of the rows in the global observability matrix
$\Bar{O}$ (ref~\eqref{eqn:obar_expand}):
\begin{align}
    O = \begin{bmatrix} (C^{1})^T\ldots(C^1 \Ax^{d-1})^T \ldots (C^{m})^T\ldots (C^m \Ax^{d-1})^T \end{bmatrix}^T. \label{eqn:o_def} 
\end{align}
Upon substituting from~\eqref{eqn:loc_obsv} in~\eqref{eqn:o_def} we have
\begin{align}
    O = \begin{bmatrix} (O^{1})^T & \ldots & (O^{m})^T\end{bmatrix}^T. \label{eqn:o}
\end{align}
Upon substituting from~\eqref{eqn:loc_out} and~\eqref{eqn:dyn_iter} in~\eqref{eqn:y_def} we obtain that
\begin{align}
   Y & = \begin{bmatrix} (C^{1})^T\ldots(C^1 \Ax^{d-1})^T \ldots (C^{m})^T\ldots (C^m \Ax^{d-1})^T \end{bmatrix}^T z(0) \nonumber \\
   & \overset{\eqref{eqn:o_def}}{=} O z(0). \label{eqn:y}
\end{align}
From the definitions of $Y$ and $O$, respectively in~\eqref{eqn:y_def} and $\eqref{eqn:o_def}$, it follows that the set of algebraic equations $\Bar{Y} = \Bar{O}z(0)$ in~\eqref{eqn:bar_y} is a rearrangement of the set of equations $Y = O z(0)$ in~\eqref{eqn:y}. So,~\eqref{eqn:bar_y} and~\eqref{eqn:y} are equivalent.
Upon substituting from~\eqref{eqn:yi_def} in~\eqref{eqn:y_def} we have
\begin{align}
    Y = \begin{bmatrix} (y^1)^T & \ldots & (y^m)^T \end{bmatrix}^T. \label{en:y_yi}
\end{align}
From~\eqref{eqn:o} and~\eqref{en:y_yi}, the set of equations in~\eqref{eqn:y} are equivalent to
\begin{align}
    y^i &= O^i z(0), ~ i = 1,\ldots,m. \label{eqn:y_dist}
\end{align}
For each agent $i$, we define the individual cost function as
\begin{align}
    F^i(x) = \frac{1}{2} \norm{O^i x - y^i}^2, \forall x\in\R^d.
\end{align}
Then, the following distributed linear regression problem solves~\eqref{eqn:y_dist}:
\begin{align}
  z(0) \in X^* = \arg \min_{x \in \R^d} \, \sum_{i = 1}^m F^i(x) = \arg \min_{x \in \R^d} \, \sum_{i = 1}^m \frac{1}{2} \norm{O^i x - y^i}^2. \label{eqn:reg}
\end{align}
Under Assumption~\ref{assump:obsv}, the matrix $\Bar{O}$ is full rank. Since $O$ is a rearrangement of the rows in $\Bar{O}$ (ref~\eqref{eqn:obar_expand} and~\eqref{eqn:o_def}), Assumption~\ref{assump:obsv} then implies that the matrix $O$ is full rank. From~\eqref{eqn:o}, then it follows that the optimization problem~\eqref{eqn:reg} has a unique solution.\\

Hence, under Assumption~\ref{assump:obsv}, the true solution of distributed linear regression problem~\eqref{eqn:reg} is the initial state $z(0)$ of the LTI system~\eqref{eqn:dyn}-\eqref{eqn:loc_out}. Estimates of the system states $z(t)$ at the subsequent sampling instants $t>0$ is obtained from~\eqref{eqn:dyn_iter}.\\

{\bf Note}: State feedback controller for LTI systems is designed by essentially solving a set of linear algebraic equation, so that the roots of the characteristic equation coincide with the desired pole locations of the closed-loop system~\cite{nagrath1999control}. Thus, the distributed state feedback controller design problem can be formulated as a distributed linear regression problem, following similar steps as above.
%%%%%%%%%%%%%%%%%%%%%%%%%%%%%%%%%%%%%%%%%%%%%%%%%%%%%%%%%%%%%%%%%%%%%%%%%%%%%%%%

\addtolength{\textheight}{-12cm}

\end{document}